\documentclass[12pt]{article}
\usepackage{graphicx}
\usepackage{epstopdf}
\usepackage{amssymb}
\usepackage{amsmath,epsfig}
\usepackage[authoryear,round]{natbib} 
\newtheorem{proposition}{Proposition}

\newtheorem{proof}{Proof}

\def\bLambda{\mbox{\boldmath $\Lambda$}}

\def\btheta{\mbox{\boldmath $\theta$}}
\def\bbeta{\mbox{\boldmath $\beta$}}

\def\bSigma{\mbox{\boldmath $\Sigma$}}
\def\bw{\mbox{\boldmath $\beta$}}

\def\bmu{\mbox{\boldmath $\mu$}}

\def\veta{\mbox{\boldmath $\eta$}}
\def\bbeta{\mbox{\boldmath $\beta$}}

\def\bz{{\bf z}}

\def\bI{{\bf I}}
\def\bU{{\bf U}}
\def\bQ{{\bf Q}}
\def\bR{\mathbb{R}}

\def\bX{{\bf X}}

\def\bw{{\bf w}}

\def\bu{{\bf u}}
\def\bx{{\bf x}}
\def\bv{{\bf v}}
\def\by{{\bf y}}
\def\nn{\nonumber}

\def\v2{\vspace{0.2in}}

\begin{document}
\setcounter{page}{1}
\baselineskip=17pt
\footskip=.3in

\title{Phase transition and higher order analysis of $L_q$ regularization under dependence}

\author{Hanwen Huang \\
    {\it Department of Epidemiology and Biostatistics}\\
    {\it University of Georgia, Athens, GA 30602}\\
    huanghw@uga.edu\\\\
Peng Zeng \\
{\it Department of Mathematics \& Statistics}\\
{\it Auburn University, Auburn, AL 36849}\\
zengpen@auburn.edu\\\\
Qinglong Yang\\
{\it School of Statistics and Mathematics}\\
{\it Zhongnan University of Economics and Law}\\
{\it Wuhan, Hubei 430073, P. R. China}\\
yangqinglong@zuel.edu.cn}
\date{}

\maketitle

\begin{abstract}
We study the problem of estimating a $k$-sparse signal $\bbeta_0\in\bR^p$ from a set of noisy observations $\by\in\bR^n$ under the model $\by=\bX\bbeta+w$, where $\bX\in\bR^{n\times p}$ is the measurement matrix the row of which is drawn from distribution $N(0,\bSigma)$. We consider the class of $L_q$-regularized least squares (LQLS) given by the formulation $\hat{\bbeta}(\lambda,q)=\text{argmin}_{\bbeta\in\bR^p}\frac{1}{2}\|\by-\bX\bbeta\|^2_2+\lambda\|\bbeta\|_q^q$,
where $\|\cdot\|_q$ $(0\le q\le 2)$ denotes the $L_q$-norm. In the setting $p,n,k\rightarrow\infty$ with fixed $k/p=\epsilon$ and $n/p=\delta$, we derive the asymptotic risk of $\hat{\bbeta}(\lambda,q)$ for arbitrary covariance matrix $\bSigma$ which generalizes the existing results for standard Gaussian design, i.e. $X_{ij}\overset{i.i.d}{\sim}N(0,1)$. We perform a higher-order analysis for LQLS in the small-error regime in which the first dominant term can be used to determine the phase transition behavior of LQLS. Our results show that the first dominant term does not depend on the covariance structure of $\bSigma$ in the cases $0\le q\textless 1$ and $1\textless q\le 2$ which indicates that the correlations among predictors only affect the phase transition curve in the case $q=1$ a.k.a. LASSO. To study the influence of the covariance structure of $\bSigma$ on the performance of LQLS in the cases $0\le q\textless 1$ and $1\textless q\le 2$, we derive the explicit formulas for the second dominant term in the expansion of the asymptotic risk in terms of small error. Extensive computational experiments confirm that our analytical predictions are consistent with numerical results.
\end{abstract}
 
\noindent%

\section{Introduction}
The goal of linear regression is to estimate the parameter vector $\bbeta_0\in\bR^p$ from a set of $n$ response variables $\by\in\bR^n$ under the model 
\begin{eqnarray}\label{lqls0}
\by=\bX\bbeta_0+\bw, 
\end{eqnarray}
where $\bX\in\bR^{n\times p}$ is a known measurement matrix and $\bw\in\bR^n$ is a noise vector. We consider the popular class of $L_q$-regularized least square methods (LQLS), given by the optimization problem
\begin{eqnarray}\label{lqls}
\hat{\bbeta}=\text{argmin}_{\bbeta\in\bR^p}\frac{1}{2}\|\by-\bX\bbeta\|^2_2+\lambda\|\bbeta\|_q^q,
\end{eqnarray}
where $\|\cdot\|_q$ $(0\le q\le 2)$ denotes the $L_q$-norm and $\lambda\ge 0$ is a fixed number. 

In many modern statistical applications with large $p$, the true $\bbeta_0$ vector is often sparse, i.e. only $k$ of the elements of $\bbeta_0$ are non-zero and the rest are zero. We would like to identify the useful predictors and also obtain good estimates of their coefficients. Among all LQLSs, the case $q=1$, or in other words $L_1$ optimization a.k.a. LASSO, is the most popular and best studied scheme for identifying a set of relevant features from a given list. LASSO enjoys attractive statistical properties \citep{10.1214/aos/1015957397,https://doi.org/10.1002/cpa.20132,10.1214/009053606000000281} and has been effectively used for simultaneously producing accurate and parsimonious models. 

However, as a convex relaxation of $L_0$-norm, the $L_1$-norm may be sub-optimal for recovering a real sparse signal because it tends to result in biased estimation by shrinking all the entries toward to zero simultaneously, and sometimes leads to over-penalization. Many researchers have shown that using the $L_q$-norm $(0\textless q\textless 1)$ to approximate the $L_0$-norm is in fact a better choice than using the $L_1$-norm because the former provides a tighter optimization relaxation to the original $L_0$ sparsity finding formulation. For example, in the noiseless settings ($w=0$), \cite{Chartrand_2008,4518502,SAAB201030} have shown that the global minima of (\ref{lqls}) for $q\textless 1$ outperforms the solution of LASSO. They have derived the Restricted Isometry Property conditions that can guarantee the accurate recovers of the sparse vector from (\ref{lqls}).

To understand the behavior of these algorithms in high dimension more reliably, a new asymptotic has been proposed recently which assumes that $p,n,k\rightarrow\infty$ with fixed $k/p=\epsilon$ and $n/p=\delta$. One of the main studies undertaken in this asymptotic framework is to find the relationship between the sparsity of the vector $\bbeta_0$ and the successful recovery of it using $L_q$-minimization algorithm (\ref{lqls}) under the noiseless setting. In other word, for different values of $\delta$, we need to know the boundary of sparsity $\epsilon$ that differentiates the success and failure of recovering by $L_q$-minimization. This boundary is known as the phase transition curve. 

\cite{https://doi.org/10.48550/arxiv.1306.3976} characterizes the phase transition curve for the case of $q=0$. \cite{5893944} analyzes the phase transition curve for $\delta\rightarrow 1$. \cite{Kabashima_2009} and \cite{NIPS2009_bdb106a0} derives the exact value of phase transition curve for any value of $0\le q\le 1$ and any value of $\delta$ using the non-rigorous replica method. In addition to the sharp characterization of the phase transition in noiseless case, \cite{7541384,7987040} also present accurate calculation of the mean squared error (MSE) in the presence of noise and compare the accuracy of estimator (\ref{lqls}) for different values of $q$ under the optimal tuning of the parameter $\lambda$. They observe that if the measurement noise $w$ is zero or small, then the global minima of (\ref{lqls}) for $q\textless 1$ (when $\lambda$ is optimally picked) outperforms the solution of LASSO with optimal $\lambda$. Furthermore, all values of $q\textless 1$ have the same performance when $w=0$. When $w$ is small, LQLS with the value of $q$ closer to 0 has a better performance. This coincides with the intuition that the performance of $L_q$- minimization is improved when $q$ decreases since it is closer to $L_0$-minimization. The analysis of \cite{7987040} is based on approximate message passing (AMP) algorithm and replica method. 

The phase transition of LQLS in the case $1\textless q\le 2$, also known as bridge regression, has been thoroughly studied in last several decades with different techniques, such as combinatorial geometry \citep{Donoho9446}, statistical dimension framework \citep{Amelunxen2013LivingOT}, Gordon's lemma \citep{8365826}, and AMP \citep{dmm,10.1214/17-AOS1651,10.1214/19-AOS1906,8618402}. It was summarized in \cite{10.1214/17-AOS1651,10.1093/imaiai/iay024,10.1093/imaiai/iaab025} that the phase transition curve is $\delta_q(\epsilon)=1$ if $2\ge q\textgreater 1$ and $\delta_q(\epsilon)=M_1(\epsilon)$ if $q=1$, where $M_1(\epsilon)$ is an increasing function with $M_1(0)=0$ and $M_1(1)=1$. Similar to \cite{7987040}, they also propose a higher-order analysis for LQLS in the small-error regime in which the first dominant term is the one that specifies the phase transition curve and the second dominant term can be used to compare the performance of different values of $q$. They observed that the actual MSE is smaller than the one predicted by the first-order term.

All the above results are based on the i.i.d. Gaussian design assumption i.e. $X_{i,j}\sim N(0,1/n)$. Our aim in this paper is to study the phase transition under arbitrary covariance dependence, i.e. $\bX$ consists of i.i.d. Gaussian rows $\bx_i\sim N(0,\bSigma)$ with covariance matrix $\bSigma\succ 0$ and $\bSigma\ne\bI_p$. The phase transition curve in the case $q=1$, i.e. LASSO, under arbitrary covariance dependence has been studied in \cite{lasso}. It was found that the LASSO phase transition curve changes with the correlation coefficients when the signed patterns of the nonzero components of $\bbeta_0$ are not symmetric. In current paper, we focus on $0\le q\textless 1$ and $1\textless q\le 2$. Toward this goal, we first derive the limiting prediction risk of LQLS estimator (\ref{lqls}) under asymptotic setting $p,n\rightarrow\infty$ with fixed $n/p=\delta$ using the non-rigorous, but widely accepted replica method from statistical physics. Then we apply the higher order analysis to the asymptotic prediction risk in terms of small noise error for a simple block diagonal covariance matrix structure. 

Here is the summary of our results: 1. The first dominant term is not influenced by the actual covariance matrix $\bSigma$ in both cases of $0\le q\textless 1$ and $1\textless q\le 2$. This is in sharp contrast to the case $q=1$ and indicates that the phase transition curve is $\delta=\epsilon$ if $0\le q\textless 1$ and $\delta=1$ if $1\textless q\le 2$ for any covariance structure; 2. The second dominant term depends on the correlation coefficient $\rho$ which is positive in the case $0\textless q\textless 1$ and negative in the case $1\textless q\le 2$. To the best of our knowledge, this is the first result to illustrate the phase transition and higher analysis for LQLS estimators under design matrices $\bX$ that have non-independent entries. 

The rest of this paper is organized as follows: In Section \ref{method}, we present the asymptotic framework of our analysis from which the distributional limit of estimator (\ref{lqls}) is derived. In Section \ref{smallq}, we study the phase transition and higher analysis for the case $0\le q\textless 1$. In Section \ref{largeq},  we study the phase transition and higher analysis for the case $1\textless q\le 2$. In Section \ref{numeric}, we show some simulation results to verify that our analytic derivations are indeed correct. The last section is devoted to the conclusion and the proofs of our main propositions are placed in Appendix.

\section{Distributional limit of LQLS estimator}\label{method}

\begin{figure}[hbtp]
	\begin{center}
		\epsfig{file=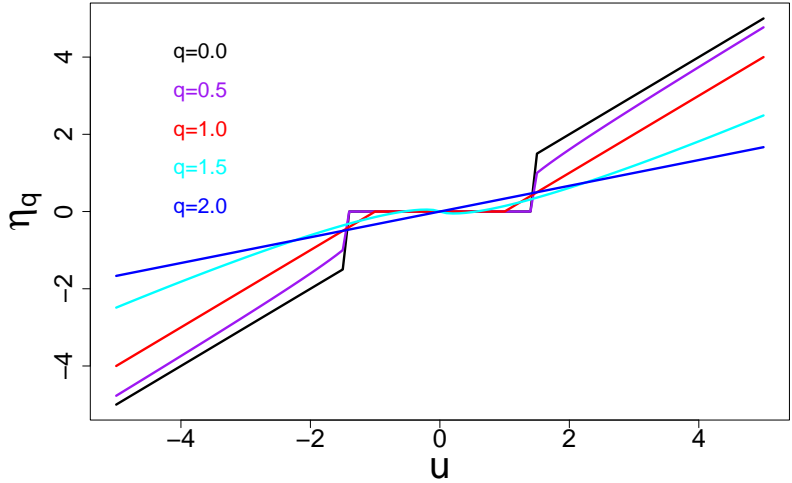,width=8.1cm,angle=-90}
	\end{center}
	\vspace{-0.5cm}
	\caption{$\eta_q(u,\alpha)$ for 5 different values of $q$. $\alpha$ is set to 1.}
	\label{figureq}
\end{figure}

The main goal of this section is to formally introduce the asymptotic setting under which we study the prediction risk of LQLS estimator (\ref{lqls}). We write a sequence of instances $\{\bbeta_0(p),\bw(p),\bSigma(p),\bX(p)\}$ indexed by $p$ which is called a converging sequence if the following conditions hold:
\begin{enumerate}
	\item $p,n\rightarrow\infty$ with $n/p\rightarrow\delta$ for some positive constant $\delta$.
	\item The empirical distribution of the entries of $\bbeta_0(p)$ converges weakly to a probability measure $p_{\beta_0}$ on $\bR$ with $\sum_{i=1}^p\beta_{0,i}(p)^2/p\rightarrow E_{p_{\beta_0}}\{\beta_0^2\}\textless\infty$.
	\item The empirical distribution of the entries of $\bw(p)$ converges weakly to a probability measure $p_{w}$ on $\bR$ with $\sum_{i=1}^nw_{i}(p)^2/n^2\rightarrow\sigma_w^2\textless\infty$.
	\item Denote $\lambda_{min}(\bSigma(p))$ and $\lambda_{max}(\bSigma(p))$ the smallest and largest eigenvalues of $\bSigma(p)$. Then $1/\lambda_{min}(\bSigma(p))=O(1)$ and $\lambda_{max}(\bSigma(p))=O(1)$.
	\item The rows of $\bX(p)$ are drawn independently from distribution $N(0,\bSigma(p)/n)$.  
	\item The sequence of functions
	\begin{eqnarray}\label{cond5}
	{\cal E}^{(p)}(a,b)&\equiv&\frac{1}{p}E\min_{\bbeta\in\bR^p}\left\{\frac{1}{2}\|\bbeta-\bbeta_0(p)-\sqrt{a}\bSigma(p)^{-1/2}\bz\|^2_{\bSigma(p)}+b\|\bbeta\|_q^q\right\}
	\end{eqnarray}
	admits a differentiable limit ${\cal E}(a,b)$ on $\bR_+\times\bR_+$ with $\text{div}{\cal E}^{(p)}(a,b)\rightarrow\text{div}{\cal E}(a,b)$, where $\|\bv\|^2_{\bSigma}=\bv^T\bSigma\bv$ and the expected value is with respect to two independent random vectors $\bz\sim N(0,\bI_{p\times p})$ and $\bbeta_0(p)$.
\end{enumerate} 
Assume we observe training data $(\bX(p),\by(p))$, where $\by(p)=\bX(p)\bbeta_0(p)+\bw(p)$. We insist on the fact that $\bbeta_0(p)$, $\bw(p)$, $\bSigma(p)$, $\bX(p)$, $\by(p)$ depend on $p$. However, we will drop this dependence most of the time to ease the reading. In order to present our main result, for any $\bv\in\bR^p$ and $\alpha\textgreater 0$, we need to introduce the generalized thresholding operation $\veta_q:\bR^p\rightarrow\bR^p$ which is defined as
\begin{eqnarray}\label{veta}
\veta_q(\bv,\alpha)=\text{argmin}_{\bbeta\in\bR^p}\left\{\frac{1}{2}\|\bbeta-\bv\|^2_{\bSigma}+\alpha\|\bbeta\|_q^q\right\}.
\end{eqnarray} 
It can be easily verified that $\veta_q(\bv,\alpha)\xrightarrow{\alpha\rightarrow 0}\bv$ and $\veta_q(\bv,\alpha)\xrightarrow{\alpha\rightarrow\infty}0$. For $\bSigma=\bI_{p\times p}$, each component of $\veta_q(\bv,\alpha)$ can be solved independently using the corresponding scalar thresholding operator $\eta_q(u,\alpha)=\text{argmin}_{\beta\in\bR}\left\{\frac{1}{2}(\beta-u)^2+\alpha|\beta|^q\right\}$ whose solution is
\begin{eqnarray}\nn
\eta_q(u,\alpha)&=&\left\{\begin{array}{ccc}uI(|u|\textgreater\sqrt{2\alpha})&if&q=0\\
\text{sign}(u)\tilde{\beta}I[|u|\textgreater C_q\alpha^{1/(2-q)}]&if&0\textless q\textless 1\\
\text{sign}(u)(|u|-\alpha)I(|u|\textgreater\alpha)&if&q=1\\
\text{sign}(u)\bar{\beta}&if&1\textless q\textless 2\\
u/(1+2\alpha)&if&q=2\end{array}\right.
\end{eqnarray} 
where $C_q=[2(1-q)]^{1/{(2-q)}}+q[2(1-q)]^{(q-1)/(2-q)}$, $\tilde{\beta}$ is the largest solution of $\beta+\frac{q\alpha}{\beta^{1-q}}=|u|$, and $\bar{\beta}$ is the solution of $\beta+q\alpha\beta^{q-1}=|u|$. Figure \ref{figureq} exhibits $\eta_q$ for different values of $q$.

For a converging sequence of instances, we can define the function
\begin{eqnarray}\label{psi}
\psi_q(\tau^2,\lambda)&=&\sigma^2_w+\lim_{p\rightarrow\infty}\frac{1}{p\delta}E\left(\|\veta_q(\bbeta_0+\tau\bSigma^{-1/2}\bz,\lambda)-\bbeta_0\|_{\bSigma}^2\right),
\end{eqnarray}
where $\bz\sim N(0,\bI_{p\times p})$ is independent of $\bbeta_0$. Notice that the function $\psi_q(\cdot,\cdot)$ depends implicitly on the law $p_{\beta_0}$.

Denote $\hat{\bbeta}(\lambda,p)$ the LQLS estimator for instance $\{\bbeta_0(p),\bw(p)$, $\bSigma(p),\bX(p)\}$ with $\lambda\textgreater 0$ based on (\ref{lqls}). Then the following proposition establishes the distributional limit of $(\hat{\bbeta}(\lambda,p),\bbeta_0)$. 
\begin{proposition}\label{prop0}
	$(\hat{\bbeta}(\lambda,p),\bbeta_0)$ converges in distribution to the random vector $(\veta_{q}(\bbeta_0+\tau_\star\bSigma^{-1/2}\bz,\alpha),\bbeta_0)$ as $p\rightarrow\infty$, where $\bz\sim N(0,\bI_{p\times p})$ is independent of $\bbeta_0\sim p_{\beta_0}$, and $\alpha,\tau_\star^2$ satisfy the following equations
	\begin{eqnarray}\label{tuning1}
	\tau_\star^2&=&\psi_q(\tau_\star^2,\alpha),\\\label{tuning2}
	\lambda&=&\alpha\left(1-\frac{1}{\delta}E\{\langle\veta_q^\prime(\bbeta_0+\tau_\star\bSigma^{-1/2}\bz,\alpha)\rangle\}\right),
	\end{eqnarray}
	where the expectation is with respect to $\bz\sim N(0,\bI_{p\times p})$ and $\bbeta_0$, $\veta_q^\prime(\cdot,\cdot)$ is the derivative of the generalized thresholding function over its first argument, and $\langle\bv\rangle\equiv\sum_{j=1}^pv_j/p$ is the average of the vector $\bv\in\bR^p$. 
\end{proposition}
The derivation of this proposition is presented in the supplementary material using the replica method which is a non-rigorous technique invented in statistical physics to study the behavior of large magnetic and disordered systems. This method has been used to analyze the accuracy of $\hat{\bbeta}(\lambda,p)$ for i.i.d. Gaussian $\bX$, i.e. $\bSigma=\bI_p$ in \cite{NIPS2009_bdb106a0}. Here we generalize it to arbitrary $\bSigma$.

The performance of LQLS estimator (\ref{lqls}) depends on the choice of the threshold parameters $\lambda$. Any fair comparison between LQLS for different values of $q$ must take this fact into account. Consider a test point $\bx_0\sim N(0,\bSigma)$ independent of the train data. For an estimator $\hat{\bbeta}$ (a function of the training data $\bX,\by$), we define its prediction risk (or simply, risk) as 
\begin{eqnarray}\nn
R_q(\hat{\bbeta}(\lambda))=\frac{1}{p}E(\bx_0^T\hat{\bbeta}(\lambda)-\bx_0^T\bbeta_0)=\frac{1}{p}E\|\hat{\bbeta}(\lambda)-\bbeta_0\|_{\bSigma}^2.
\end{eqnarray}
For $\bSigma=\bI_p$, $R_q(\hat{\bbeta}(\lambda))$ is equivalent to the MSE defined in \cite{7541384}. Hence, it is important to pick $\lambda$ optimally by minimizing the prediction risk, i.e. $\lambda_\star=\text{argmin}_\lambda R_q(\hat{\bbeta}(\lambda))$. This is equivalent to consider the following optimal thresholding policy 
\begin{eqnarray}\label{psi}
\alpha_\star(\tau)=\text{argmin}_\alpha\lim_{p\rightarrow\infty}\frac{1}{p}E\left(\|\veta_q(\bbeta_0+\tau\bSigma^{-1/2}\bz,\alpha)-\bbeta_0\|_{\bSigma}^2\right)
\end{eqnarray}
since $R_q(\hat{\bbeta}(\lambda_\star))=\delta[\psi(\tau_\star^2,\alpha_\star(\tau_\star))-\sigma^2_w]$.
This enables us to focus the analysis on a single equation rather than two equations (\ref{tuning1}) and (\ref{tuning2}). The results we will present in the next two sections are mainly based on investigating the solution of the following fixed point equation
\begin{eqnarray}\label{fixed}
\tau_\star^2&=&\psi_q(\tau_\star^2,\alpha_\star(\tau_\star)).
\end{eqnarray}
We consider the cases $0\le q\textless 1$ and $1\textless q\le 2$ separately.

\section{Asymptotic performance for $0\le q\textless 1$}\label{smallq}
When $0\le q\textless 1$, the $L_q$ regularization has the nice property of creating coefficient sparsity because the generalized thresholding operation $\veta_q(\bv,\lambda)$ defined in (\ref{veta}) maps small values of $\bv$ to zero. Also the smaller values of $q$ correspond to a preference for increasingly sparse solutions. However, the penalty term is not convex for $q\textless 1$ which makes solving the optimization problem much more difficult because the local minima may not be unique.
\subsection{Noiseless Settings}\label{noislesspl}
This section discusses our main results in the noiseless setting $\sigma_w^2=0$. Since there is no measurement noise, $\psi_q(\tau^2,\alpha_\star(\tau))\rightarrow 0$ as $\tau^2\rightarrow 0$, thus $\tau^2=0$ is a fixed point of (\ref{fixed}). If $\tau^2=0$ is also a stable fixed point of (\ref{fixed}), i.e. there exists $\tau_i\textgreater 0$ such that for every $0\textless\tau\textless\tau_i$, $\psi_q(\tau^2,\alpha_\star(\tau))\textless\tau^2$, then we say LQLS has successfully recovered the sparse solution of $\by=\bX\bbeta_0$. Otherwise, we say LQLS has failed.

The shape of $\psi_q(\tau^2,\alpha_\star(\tau))$ and its fixed points depend on the distribution $p_{\beta_0}$. In this work we focus on $p_{\beta_0}\in{\cal F}_\epsilon$, where ${\cal F}_\epsilon$ denotes the set of distributions whose mass at zero is greater than or equal to $1-\epsilon$. In other words, $\beta_0\sim p_{\beta_0}$ implies that $P(\beta_0\ne 0)\le\epsilon$. This class of distributions is considered as a good model for exactly sparse signals and has been studied in many other papers. The following proposition identifies the conditions under which zero is a stable fixed point under noiseless setting.
\begin{proposition}\label{prop1}
	Let $p_{\beta_0}$ be an arbitrary distribution in ${\cal F}_\epsilon$. For any $0\le q\textless 1$, zero is the stable fixed point of $\psi_q(\tau^2,\alpha_\star(\tau))$ under noiseless settings if and only if $\delta\textgreater\epsilon$. Thus the phase transition curve is $\delta=\epsilon$.
\end{proposition}
The proof of this result can be found in Section \ref{prop1p}. This result extends the conclusion of Theorem 4 in \cite{7987040} from $\bSigma=\bI_{p\times p}$ to general $\bSigma$. There are three main features of this proposition that we would like to emphasize. Firstly, the actual distribution that is picked from ${\cal F}_\epsilon$ does not have any impact on the behavior of the fixed point at zero. Secondly, this proposition is universal in the sense that it does not depend on the actual structure of covariance matrix $\bSigma$. Thirdly, the number of the measurements $\delta$ that is required for the stability of this fixed point is the same as the sparsity level $\epsilon$. As long as $\delta\textgreater\epsilon$, zero is a stable fixed point and LQLS recovers $\bbeta_0$ accurately for every $0\le q\textless 1$. If we are concerned with the noiseless settings, all $L_q$-minimization algorithms are the same. 

Note that since $L_q$ is not convex for $0\le q\textless 1$, $\psi_q(\tau^2,\alpha_\star(\tau))$ may have additional stable fixed points than $\tau^2=0$. The conditions under which this case happens depend on $q$, $p_{\beta_0}$,  and the actual structure of $\bSigma$. We will provide a numeric study on this type of phase transition in Section \ref{num1} for $q=0$ and a block diagonal $\bSigma$.

As pointed out in \cite{7987040}, this result seems to be counter-intuitive. For instance, we expect to see that the performance difference between $q=0$ and $q=0.9$ is bigger than the performance difference between $q=0$ and $q=0.1$. We also expect to see the impact of covariance structure of $\bSigma$ on the performance of LQLS estimators. However, according to the phase transition analysis of Proposition \ref{prop1}, the performance of all LQLS are the same under the noiseless settings. This naturally leads to the following question: to what extent are the phase transition analysis applicable if the response variables include small amount of errors in practice, i.e. the variance $\sigma_w^2$ of the error $\bw$ is assumed to be small. In the next section, we will clarify this surprising phenomenon by performing high order analysis in terms of $\sigma_w$. 

\subsection{Noisy Settings}\label{noisypl}
In this section, we assume that $\sigma_w^2\textgreater 0$. This implies that zero is no longer a fixed point of $\psi_q(\tau^2,\alpha_\star(\tau))$ defined in (\ref{fixed}) and the reconstruction error of LQLS is greater than zero for all values of $q$. Denote $\tau_l$ the lowest fixed point of $\psi_q(\tau^2,\alpha_\star(\tau))$. Our first result is concerned with the first dominant term for small amount of noise. 
\begin{proposition}\label{prop2}
	If $\epsilon\textless\delta$, then there exists $\tau_s^2$ such that for every $\sigma_w^2\textless\tau_s^2$, $\lim_{\sigma_w^2\rightarrow 0}\frac{\tau_l^2}{\sigma_w^2}=\frac{\delta}{\delta-\epsilon}$.
\end{proposition}
As discussed in Section \ref{noislesspl}, the first dominant term determines the phase transition which is the same for all $q$ and $\bSigma$. In order to see the discrepancy between different values of $q$ and different structures of $\bSigma$, we have to do high-order analysis and explore how $\tau_l^2-\frac{\delta\sigma_w^2}{\delta-\epsilon}$ behaves for small values of $\sigma_w^2$. We consider the simplest block diagonal matrix $\bSigma$ with two-dimensional block $\left(\begin{array}{cc}1&\rho\\\rho&1\end{array}\right)$ where $0\le\rho\textless 1$. For this choice of the covariance matrix, it can be easily verified that the function (\ref{cond5}) defined in Condition 5 has a differentiable limit. 

Consider $p_{\bbeta_0}\in{\cal F}_\epsilon$. Let $p_{\beta_0}=(1-\epsilon)\delta_0+\epsilon G$ and $U$ is a random variable with distribution $G$, where $\delta_0$ denotes a point mass at zero and $G$ denotes the distribution of nonzero elements of $\bbeta_{0}$. Let us discuss the result for $q=0$ and $0\textless q\textless 1$ separately in the following two propositions.
\begin{proposition}\label{prop3}
	Suppose $E|U|^2\textless\infty$ and $P(|U|\textgreater\mu)=1$, where $\mu=\sup_v\{v:P(|U|\textgreater v)=1\}\textgreater 0$. Then for $q=0$ and $\epsilon\textless\delta$, 
	\begin{eqnarray}\label{p0sl}
	\lim_{\sigma_w\rightarrow 0}\tau_l^2=\frac{\delta}{\delta-\epsilon}\sigma_w^2+o\left(\phi\left(\tilde{\mu}\sqrt{\frac{\delta-\epsilon}{\delta}}\sigma_w^{-1}\right)\right), 
	\end{eqnarray} 
	where $\phi$ denotes the standard normal density function and $\tilde{\mu}$ is any constant that is smaller than $(2-\rho-\sqrt{1-\rho^2})\sqrt{1-\rho^2}\mu/2$.
\end{proposition}
\begin{proposition}\label{prop4}
	Suppose $E|U|^2\textless\infty$ and $P(|U|\textgreater\mu)=1$, where $\mu=\sup_v\{v:P(|U|\textgreater v)=1\}\textgreater 0$. Then for $0\textless q\textless 1$ and $\epsilon\textless\delta$, 
	\begin{eqnarray}\label{p1sl}
	\lim_{\sigma_w\rightarrow 0}\frac{\tau_l^2-\frac{\delta}{\delta-\epsilon}\sigma_w^2}{\sigma_w^{4-2q}(\log\frac{1}{\sigma_w})^{2-q}}&=&\frac{\epsilon C_q^{4-2q}q^2(\epsilon D_0+(1-\epsilon)C_0)\delta^{2-q}}{(4-4q)^{2-q}(\delta-\epsilon)^{3-q}},
	\end{eqnarray} 
	where $C_q=[2(1-q)]^{1/(2-q)}+q[2(1-q)]^{(q-1)/(2-q)}$,  $D_0=\left[E\{|U|^{2q-2}\}-\rho(E\{|U|^{q-1}sign(U)\})^2\right]/(1-\rho^2)$, and $C_0=E\{|U|^{2q-2}\}$.
\end{proposition}
It can be easily verified that when $\rho=0$, results from Proposition \ref{prop3} and Proposition \ref{prop4} are the same as the results from Theorem 8 and Theorem 9 in \cite{7987040}. Comparing (\ref{p0sl}) and (\ref{p1sl}), we conclude that the second dominant term for $q=0$ decays exponentially faster than the polynomial rate for $0\textless q\textless 1$ in low noise regime. Since the second dominant term in (\ref{p1sl}) is positive, optimally tuned $L_0$ regularization will outperform optimally tuned $L_q$ regularization for $0\textless q\textless 1$ in this regime. Another interesting feature of Proposition \ref{prop4} is that the second dominant term is proportional to $\sigma_w^{4-2q}$, thus if $q_1\textless q_2$, $LQLS$ for $q_1$ outperforms $LQLS$ for $q_2$ for small enough $\sigma_w^2$. Moreover, the second dominant term increases with $\epsilon$ and decreases with $\delta$. Examining the impact of $\rho$ based on (\ref{p1sl}), we conclude that if the distribution of U is symmetric about zero, the second dominant term increases with $\rho$ since $D_0$ is proportional to $1/(1-\rho^2)$ and $E\{|U|^{q-1}sign(U)\}=0$. All these observations are consistent with the numerical studies shown in Section \ref{numeric}.

\section{Asymptotic Performance for $1\textless q\le 2$}\label{largeq}
The difference between the regularization for $1\textless q\le 2$ and the regularization for $0\le q\textless 1$ is that the former is convex and its optimization problem has a unique global minimum. However, since its penalty is differentiable everywhere, it never leads to a coefficient been zero rather only minimizes it when $\lambda\ne 0$. Our first result is concerned with the first dominant term of the optimally tuned asymptotic prediction risk defined as $R_q(\tau_\star,\sigma_w)=\min_\lambda\frac{1}{p}\|\hat{\bbeta}(\lambda)-\bbeta_0\|_{\bSigma}^2$ which is related to the phase transition curve for successful recovery. 
\begin{proposition}\label{prop5}
	For $\bbeta_0\in{\cal F}_\epsilon$, suppose $E|U|^2\le\infty$ and $P(|U|\textgreater\mu)=1$, where $\mu=\sup_v\{v:P(|U|\textgreater v)=1\}\textgreater 0$. Then for $1\textless q\le 2$, we have 
	\begin{eqnarray}\nn
	\lim_{\sigma_w^2\rightarrow 0}R_q(\tau_\star,\sigma_w)\left\{\begin{array}{ccc}\textgreater 0&if&\delta\textless 1\\\rightarrow 0&if&\delta\textgreater 1
	\end{array}.\right.
	\end{eqnarray}	
\end{proposition}
Therefore, assume the error $\sigma_w^2$ equals zero, for $1\textless q\le 2$ and any $\gamma\textgreater 0$, if $\delta\textgreater 1+\gamma$, then (\ref{lqls}) succeeds in recovering $\bbeta$, while if $\delta\textless 1-\gamma$, (\ref{lqls}) fails. The phase transition curve $\delta=1$ is independent of $\epsilon$, $q$, and the covariance structure of $\bSigma$. As pointed out by \cite{7987040}, phase transition analysis based on the first dominant term may lead to misleading conclusions in any practical setting where there is always error in the response variables. To study the impact of them, we need to conduct higher-order analysis of LQLS in the small-error regime in the expansion of prediction risk $R_q(\tau_\star,\sigma_w)$.

Our second result is concerned with the second dominant term of the optimally tuned MSE of LQLS when the number of response variables is larger than the number of predictors $p$, i.e. $\delta\textgreater 1$. The key is to characterize the convergence rate for $R_q(\tau_\star,\sigma_w)$ as $\sigma_w\rightarrow 0$. 
\begin{proposition}\label{prop6}
	Denote $Z\sim N(0,1)$ and $U$ a random variable whose distribution is specified by the nonzero elements of $\beta_0$. Suppose $E|U|^2\textless\infty$ and $P(|U|\textgreater\mu)=1$, where $\mu=\sup_v\{v:P(|U|\textgreater v)=1\}\textgreater 0$. Then for $1\textless q\le 2$ and $\delta\textgreater 1$, we have 
	\begin{eqnarray}\nn
	&&R_q(\tau_\star,\sigma_w)\\\label{thanone}
	&=&\frac{\delta\sigma_w^2}{\delta-1}-\frac{\delta^{q+1}}{(\delta-1)^{q+1}}\frac{(1-\epsilon)^2}{\epsilon}\frac{[E(|Z|^{q})]^2\sigma_w^{2q}}{(1-\rho^2)^{q-1}\{E(|U|^{2q-2})-\epsilon\rho[E(|U|^{q-1}\text{sign}(U))]^2\}}.
	\end{eqnarray}	
\end{proposition}
The first term $\frac{\delta\sigma_w^2}{\delta-1}$ shows a notion of phase transition. For as $\sigma_w\rightarrow 0$, $R_q(\tau_\star,\sigma_w)=O(\sigma_w^2)$, and will go to zero. For $\rho=0$, the second dominant term is the same as the one derived in Theorem 3.1 in \cite{10.1214/17-AOS1651}. The important facts of the second dominant term include: (1) it is negative; (2) its order is $\sigma_w^{2q}$; (3) its magnitude decreases with $\epsilon$ for fixed $q$ and $\rho$; (4) its magnitude decreases with $q$ for fixed $\epsilon$ and $\rho$; (5) its magnitude increases with $\rho$ for fixed $q$ and $\epsilon$. Facts (3) and (4) have been studied thoroughly in \cite{10.1214/17-AOS1651}. Our numerical studies in Section 5 focus on the verification of fact (5).

\section{Numerical results}\label{numeric}
This section performs several numerical studies to evaluate the $R_q(\tau_\star,\sigma_w)$ of LQLS as a function of $\sigma_w$ for several different values of $q$ and $\rho$. Here we consider a block diagonal matrix $\bSigma$ with two-dimensional block $\left(\begin{array}{cc}1&\rho\\\rho&1\end{array}\right)$. The equation can be decomposed into 2-dimensional blocks. Sections \ref{num1} and \ref{num2} study the performance of LQLS with $0\le q\textless 1$ and $1\textless q\le 2$ respectively. 
\subsection{$0\le q\textless 1$}\label{num1}
\begin{figure}[hbtp]
	\begin{center}
		\epsfig{file=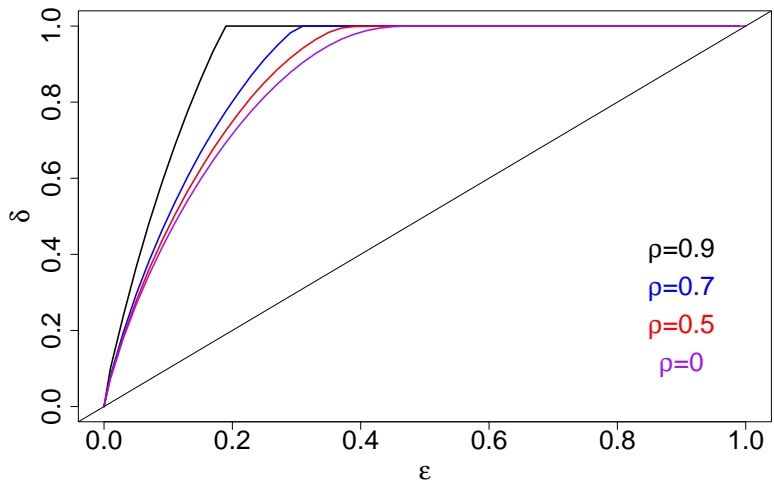,width=8.1cm,angle=-90}
	\end{center}
	\vspace{-0.5cm}
	\caption{Comparison of the phase transition of the $L_0$ regularized method. The phase transition exhibited is the value of $\delta$ at which the number of stable fixed points changes from one to more than one for at least some prior $P_{\beta_0}\in{\cal F}$.}
	\label{figurep00phase}
\end{figure}

\begin{figure}[hbtp]
	\begin{center}
		\epsfig{file=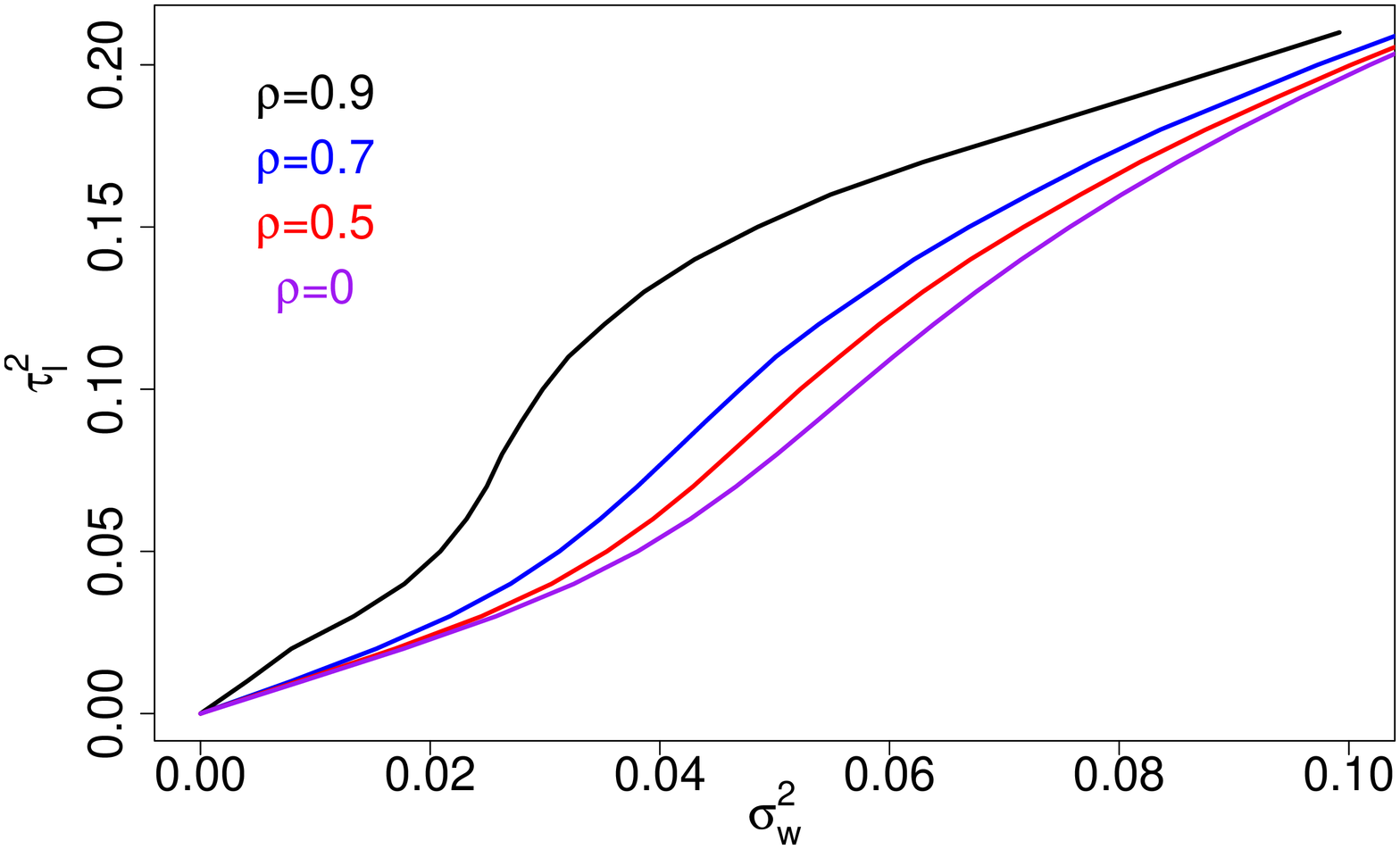,width=8.1cm,angle=-90}
	\end{center}
	\vspace{-0.5cm}
	\caption{The curve of $\tau_l^2$ as a function of $\sigma_w^2$ for $\rho\in[0,0.5.0.7.0.9]$. Here $\delta=0.9$, $\epsilon=0.1$, and the non-zero elements of $\beta_0$ are i.i.d. $\pm 1$ with probability 0.5.}
	\label{figurep00}
\end{figure}

\begin{figure}[hbtp]
	\begin{center}
		\epsfig{file=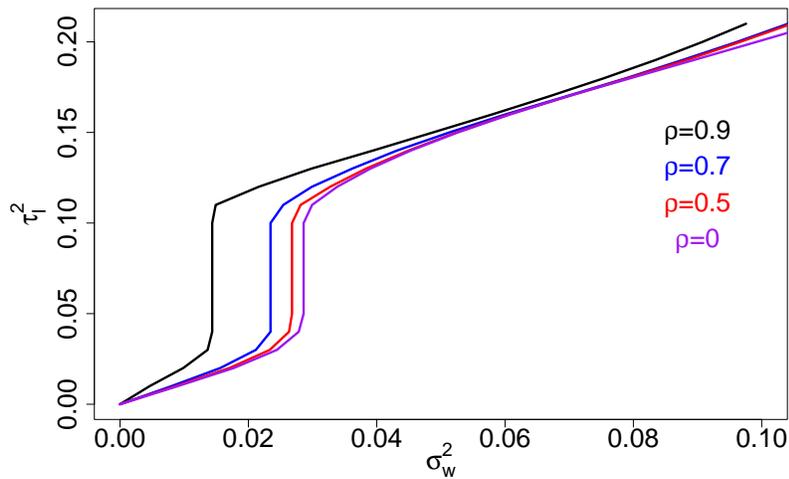,width=8.1cm,angle=-90}
	\end{center}
	\caption{The curve of $\tau_l^2$ as a function of $\sigma_w^2$ for $\rho\in[0,0.5.0.7.0.9]$. Here $\delta=0.1$, $\epsilon=0.01$, and the non-zero elements of $\beta_0$ are i.i.d. $\pm 1$ with probability 0.5.}
	\label{figurep00d1e01}
\end{figure}

We first study the phase transition result for $q=0$, i.e. $\sigma_w=0$, in Figure \ref{figurep00phase} to identify conditions under which $\psi_q(\tau^2,\alpha_\star(\tau))$ has additional stable fixed points than $\tau^2=0$. The following quantity plays an important role in our analysis
\begin{eqnarray}\label{m0}
M_{0}(\epsilon,\rho)=\sup_{\mu_1,\mu_2}\inf_\alpha\left\{\frac{\epsilon^2}{2}E\|\hat{\beta}_1-\bbeta_0^1\|^2_{\bSigma}+\epsilon(1-\epsilon)E\|\hat{\beta}_2-\bbeta_0^2\|^2_{\bSigma}+\frac{(1-\epsilon)^2}{2}E\|\hat{\beta}_3\|^2_{\bSigma}\right\},
\end{eqnarray} 
where
\begin{eqnarray}\nn
\hat{\beta}_1=\veta_0(\bbeta_0^1+\bSigma^{-1/2}\bz,\alpha),~\hat{\beta}_2=\veta_0(\bbeta_0^2+\bSigma^{-1/2}\bz,\alpha),~\hat{\beta}_3=\veta_0(\bSigma^{-1/2}\bz,\alpha),
\end{eqnarray} 
where 
\begin{eqnarray}\nn
\bSigma=\left(\begin{array}{cc}1&\rho\\\rho&1\end{array}\right), ~\bbeta_0^1=\left(\begin{array}{c}\mu_1\\\mu_2\end{array}\right),~\bbeta_0^2=\left(\begin{array}{c}\mu_1\\0\end{array}\right).
\end{eqnarray} 
Here we consider 3 different situations according to the nonzero components of 2-dimensional vector $\bbeta_0$. The first term in (\ref{m0}) is for situation that both components are nonzero. The second term is for situation that one is nonzero and the other is zero. The third term is for situation that both are zero. From (\ref{psi}) we obtain that if $\delta\textgreater M_{0}(\epsilon,\rho)$, $\psi_0(\tau^2,\lambda_\star(\tau))$ has a unique stable fixed point at zero. On the other hand, if $\delta\textless M_{0}(\epsilon,\rho)$, there exists $P_{\beta_0}\in{\cal F_\epsilon}$ for which $\psi_0(\tau^2,\lambda_\star(\tau))$ has more than one stable fixed point. Figure \ref{figurep00phase} compares the phase transition curves for different values of $\rho$. For larger $\rho$, we need larger $\delta$ to achieve unique fixed point for $\psi_0(\tau^2,\lambda_\star(\tau))$. Figures \ref{figurep00} and \ref{figurep00d1e01} exhibit the lowest fixed point $\tau_l^2$ of $\psi_0(\tau^2,\lambda_\star(\tau))$ as a function of $\sigma_w^2$ for different values of $\rho$. Figure \ref{figurep00} is for the case of $\delta\textgreater M_0(\epsilon,\rho)$ while Figure \ref{figurep00d1e01} is for the case of $\delta\textless M_0(\epsilon,\rho)$. Note that the performance of $L_0$ regularized method for smaller $\rho$ is better than its performance for larger $\rho$. Moreover, $\tau_l^2$ is a continuous function of $\sigma_w^2$ if $\psi_0(\tau^2,\lambda_\star(\tau))$ has a unique stable fixed point as shown in Figure \ref{figurep00}. On the other hand, if $\psi_0(\tau^2,\lambda_\star(\tau))$ has more than one stable fixed points, the function is discontinuous as shown in Figure \ref{figurep00d1e01}. 

\begin{figure}[hbtp]
	\begin{center}
		\epsfig{file=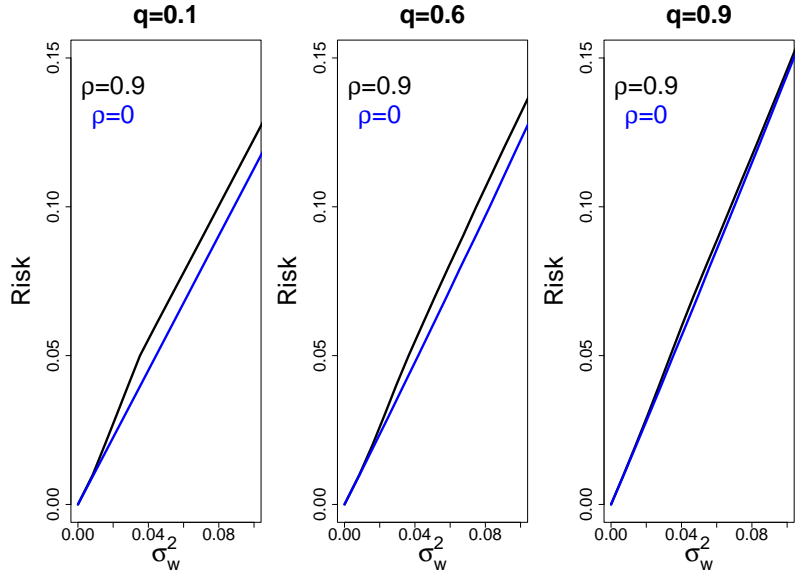,width=8.1cm,angle=-90}
	\end{center}
	\caption{The curve of MSE as a function of $\sigma_w^2$ for $\rho\in[0,0.9]$ under three different values of $q$. Here $\delta=0.9$, $\epsilon=0.1$, and the non-zero elements of $\beta_0$ are i.i.d. $\pm 3$ with probability 0.5.}
	\label{figured09e01rho}
\end{figure}

\begin{figure}[hbtp]
	\begin{center}
		\epsfig{file=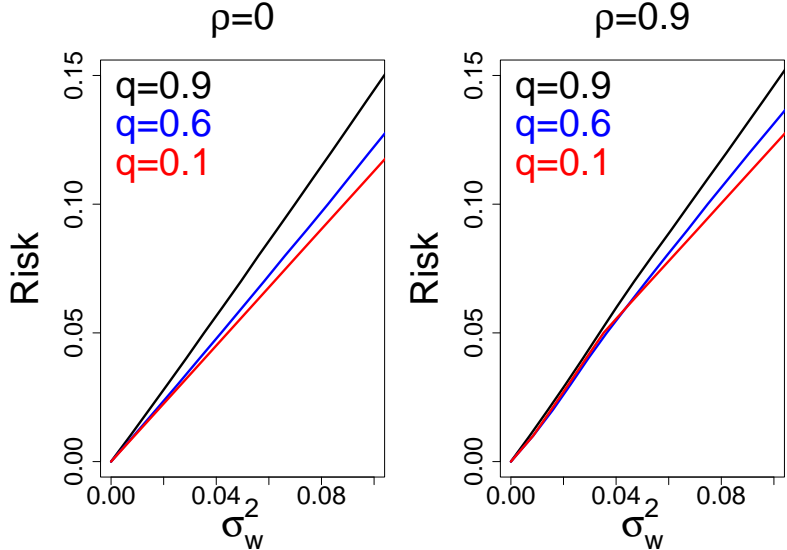,width=8.1cm,angle=-90}
	\end{center}
	\caption{The curve of MSE as a function of $\sigma_w^2$ for $q\in[0.1,0.6,0.9]$ under two different values of $\rho$. Here $\delta=0.9$, $\epsilon=0.1$, and the non-zero elements of $\beta_0$ are i.i.d. $\pm 3$ with probability 0.5.}
	\label{figured09e01q}
\end{figure}

Figure \ref{figured09e01rho} illustrates the dependence of the risk function on the correlation coefficient $\rho$ for $0\textless q\textless 1$ under different choices of $q$. It is shown clearly that the risk increases with $\rho$ for fixed $q$ which is consistent with the analytical result (\ref{p1sl}) derived in Proposition \ref{prop4} based on high order analysis. Moreover, the magnitude of changes due to $\rho$ is larger for small $q$ than that for large $q$. In particular, as $q$ close to 1, the curve for $\rho=0$ and the curve for $\rho=0.9$ in the right panel of Figure \ref{figured09e01rho} are almost indistinguishable. This is consistent with the conclusion drawn in Theorems 3.3 of \cite{10.1214/17-AOS1651} which states that the first order term is sufficient to describe the performance of LASSO (q=1) in the small-error regime because the second dominant term is exponentially small. Figure \ref{figured09e01q} depicts the impact of $q$ on risk for fixed $\rho$. It is shown that $LQLS$ for $q_1$ outperforms $LQLS$ for $q_2$ when $q_1\textless q_2$. This is consistent with the conclusion in \cite{10.1214/17-AOS1651} and also confirms our observation in Proposition \ref{prop4}. Moreover, it is shown from Figure \ref{figured09e01q} that the discrepancies caused by different values of $q$ decreases with $\rho$. 

\subsection{$1\textless q\le 2$}\label{num2}
In this section, we check the approximation accuracy of the first and second order expansions of the risk over a reasonable range of $\sigma_w$ in the case $1\textless q\le 2$. The dependence of risk on parameters $\delta$, $\epsilon$, and $q$ has been studied extensively in \cite{10.1214/17-AOS1651}. Here we focus on the dependence of risk on the correlation coefficient $\rho$. Throughout this section, we set the distribution of $U$ to $f(b)=0.5\delta_1(b)+0.5\delta_{-1}(b)$. 

\begin{figure}[hbtp]
	\begin{center}
		\epsfig{file=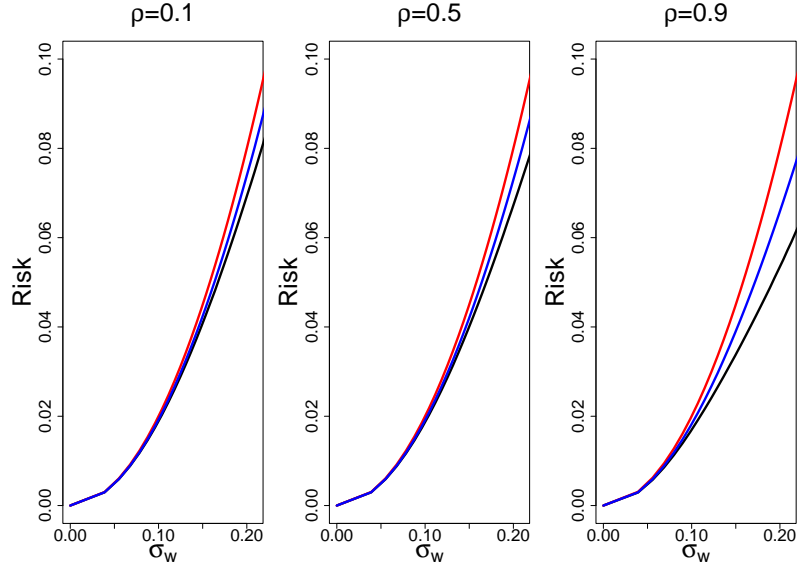,width=8.1cm,angle=-90}
	\end{center}
	\caption{Plots of actual risk and its approximation for $q=1.5,\delta=2$, and $\epsilon=0.7$. The true value, first order approximation, and second order approximation are denoted by black, red, and blue curves respectively.}
	\label{figured2p15}
\end{figure}

\begin{figure}[hbtp]
	\begin{center}
		\epsfig{file=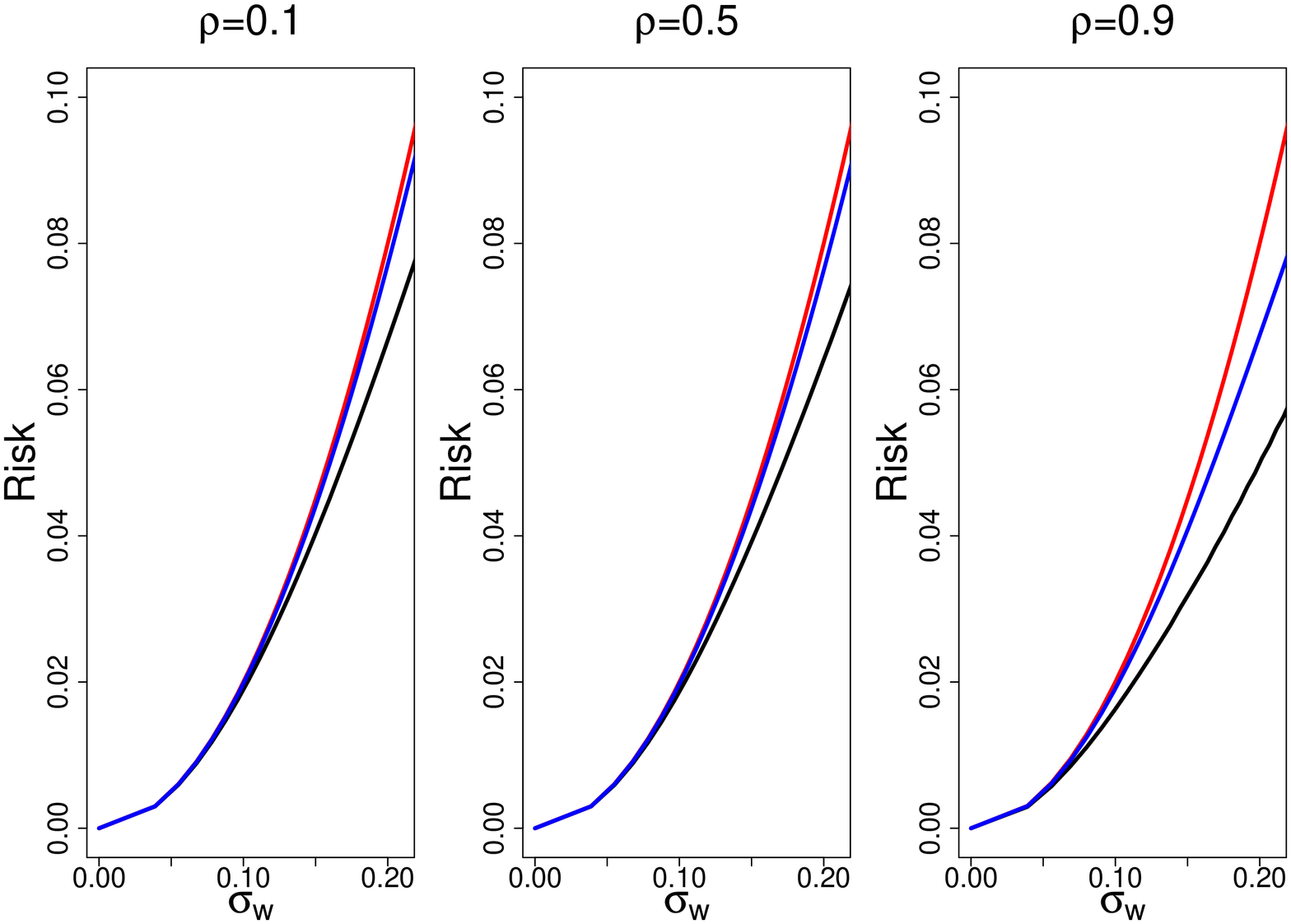,width=8.1cm,angle=-90}
	\end{center}
	\caption{Plots of actual risk and its approximation for $q=1.9,\delta=2$, and $\epsilon=0.7$. The true value, first order approximation, and second order approximation are denoted by black, red, and blue curves respectively.}
	\label{figured2p19}
\end{figure}

Figures \ref{figured2p15} and \ref{figured2p19} compare the true value, first order approximation, and second order approximation of risk under different choices of $\rho$ for $q=1.5$ and $q=1.9$ respectively. Our numerical results in Figures \ref{figured2p15} and \ref{figured2p19} show that both the first order and second order expansions present a good approximation when $\sigma_w$ is small enough. As we increase $\sigma_w$, both first order and second order expansions are larger than the true values and the second order approximation is more accurate than the first order approximation. The discrepancies between the approximation and the true values increase with $\rho$ as demonstrated in these Figures. Therefore, the expansion approximation is less accurate under large correlation than under small correlation. Consistence with the analytical form of the second-order term in (\ref{thanone}), we conclude from the numerical studies that the second order approximation outperforms the first order, however, it may not be sufficient if the following conditions hold: (i) $\delta$ is close to 1; (ii) $\epsilon$ is small, (iii) $q$ is close to 1, (iv) $\rho$ is close to 1. 

\begin{figure}[hbtp]
	\begin{center}
		\epsfig{file=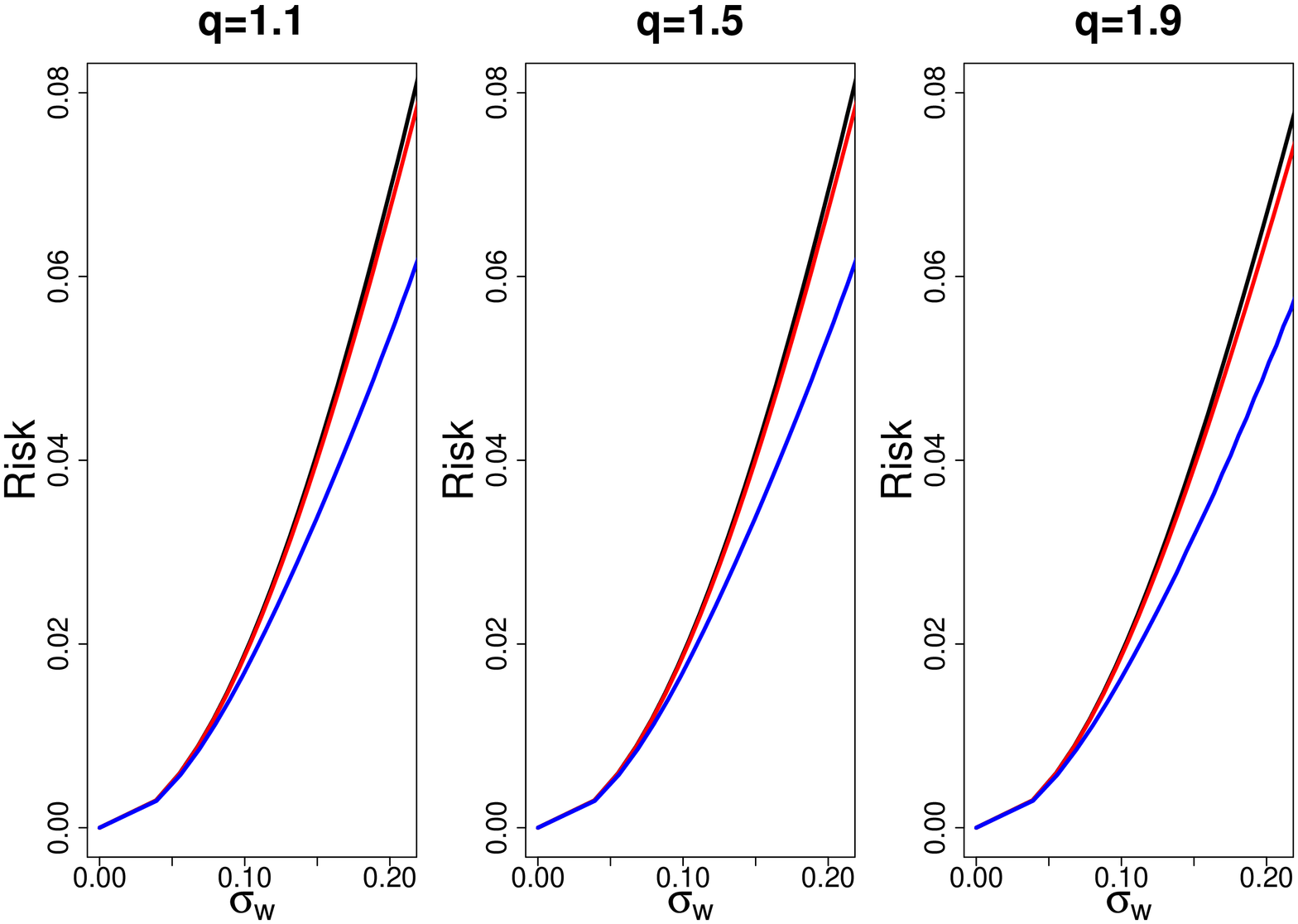,width=9.1cm,angle=-90}
	\end{center}
	\caption{Plots of actual risk for $q=1.5,\delta=2$, and $\epsilon=0.7$. The results for $\rho=0.1$, $\rho=0.5$, and $\rho=0.9$ are denoted by black, red, and blue curves respectively.}
	\label{figured2thanone}
\end{figure}

Figure \ref{figured2thanone} exhibits the impact of $\rho$ on the true value of risk for different values of $q$. As is clear in this figure, the risk decreases with both $\rho$ and $q$ which is in agreement with the second order approximation in (\ref{thanone}).  

\section{Discussion}\label{dis}

$L_q$-regularized least square is one of the most popular schemes for recovering a high-dimensional sparse vector from low-dimensional measurements. This paper focuses on asymptotic behavior of LQLS estimators in the framework in which it is assumed that $p,n,k\rightarrow\infty$ while $n/p\rightarrow\delta$ and $k/p\rightarrow\epsilon$ are fixed. We first derive the distributional limit of LQLS estimators for nonstandard Gaussian design models where the row of design matrix $\bX$ are drawn independently from distribution $N(0,\bSigma)$ with arbitrary $\bSigma$. Then we obtain an explicit characterization of the asymptotic risk by deriving its higher order expansion in the small-error regime. Our analysis is performed on both the convex case $1\textless q\le 2$ and non-convex case $0\le q\textless 1$. We conclude that the first order term does not depend on the covariance structure of $\bSigma$ and thus the phase transition curves are the same as the in the case of standard Gaussian design. This is different from the case $q=1$ in which the correlation can change the phase transition boundary if the signed pattern of signal is not symmetric. The second order term depends on the correlation and the explicit formulas are derived in both cases based on a simple two-dimensional block diagonal covariance model.  

Note that part of our analysis is based on the replica method which are not fully rigorous yet. So far the rigorous work in this area mainly focuses on i.i.d. randomness. For example, the mean square error of LQLS for the case of $1\le q\le 2$ is characterized using Gordon's lemma in \cite{8365826} and AMP in \cite{10.1214/17-AOS1651}. For nonstandard Gaussian design models, $\bSigma\ne I_{p\times p}$ the rigorous analysis for the case $q=1$ have been conducted recently \cite{wei} using Gordon’s comparison inequality and in \cite{lasso} using AMP. Our next step is to derive the rigorous results for Proposition \ref{prop0} for the case $1\textless q\le 2$. The case $0\textless q\le 1$ is more challenging because the penalty term is not convex and even in the i.i.d. case its rigorous results have not been established yet.    

Our conclusion in this paper is based on higher-order analysis in small error region. It may not be applicable if the noisy error is large enough. As shown in \cite{7987040}, for large values of noise, the LASSO outperforms $L_0$-minimization and if $q_1\textgreater q_2$, then optimal $L_{q_1}$-minimization outperforms optimal $L_{q_2}$-minimization. Proposition 3 of \cite{7987040} provides theoretical justification on this high-noise phenomenon. Another direction of our future research is to study the high-noise phenomenon of LQLS for nonstandard Gaussian design models.

\section{Appendix}
\label{sec1}
\setcounter{equation}{0}
\def\theequation{A\arabic{equation}}

This appendix outlines the proofs of the Propositions in the main text.  
 
\subsection{Proof of Proposition 1}\label{prop0p}
 
\begin{proof}
	We limit ourselves to the main steps of replica calculations leading to Propositions 1. For a general introduction to the method and its motivation, we refer to \cite{mezard1987spin,mezard2009information}.
	
	We consider $L_q$-regularized least squares estimators of the form \footnote{There is a subtle discrepancy between this definition and (\ref{lqls}) in the main text, i.e. $\bX$ and $\by$ are scaled by $1/\sqrt{n}$.}
	\begin{eqnarray}\label{lse}
	\hat{\btheta}=\text{argmin}_{\btheta\in\bR^p}\left\{\frac{1}{2n}\|\by-\bX\btheta\|_2^2+\lambda\|\btheta\|_q^q\right\}.
	\end{eqnarray}
	Define $g: \bR\times \bR\rightarrow \bR$ a continuous function strictly convex in its first argument. Then its Lagrange dual $\tilde{g}(x,y)\equiv\text{max}_{\mu\in\bR}\{\mu x-g(\mu,y)\}$ is also a continuous function convex in its first argument. We start from estimating the following moment generating function also called partition function
	\begin{eqnarray}\label{partition}
	Z_p(\beta,s)=\int\exp\left\{-\frac{\beta}{2n}\sum_{i=1}^n(y_i-\bx_i^T\btheta)^2-\beta\lambda\|\btheta\|_q^q-\beta s\sum_{j=1}^p[g(\mu_j,\beta_{0,j})-\mu_j\theta_j]\right\}d\btheta d\bmu, 
	\end{eqnarray}
	where $s\textgreater 0$, $(y_i,\bx_i)$ are i.i.d. pairs distributed as model (\ref{lqls}) in the main text, and $\beta\textgreater 0$ is the inverse temperature parameter.  The free energy is the low temperature limit
	\begin{eqnarray}
	{\cal F}(s)=-\lim_{p\rightarrow\infty}\lim_{\beta\rightarrow\infty}\frac{1}{p\beta}\log Z_p(\beta,s),
	\end{eqnarray}
	which is assumed to converge to 
	\begin{eqnarray}\label{rate}
	{\cal F}(s)=-\lim_{p\rightarrow\infty}\lim_{\beta\rightarrow\infty}\frac{1}{p\beta}E\log Z_p(\beta,s),
	\end{eqnarray}
	where the expectation is with respect to the distribution of $(y_1,\bx_1),\cdots,(y_n,\bx_n)$. Using Laplace method in the integral (\ref{partition}), we have
	\begin{eqnarray}
	{\cal F}(s)=\lim_{p\rightarrow\infty}\frac{1}{p}\min_{\btheta,\bmu\in\bR^p}\left\{\frac{1}{2n}\sum_{i=1}^n(y_i-\bx^T_i\btheta)^2+\lambda\|\btheta\|_q^q+s\sum_{j=1}^p[g(\mu_j,\beta_{0,j})-\mu_j\theta_j]\right\}. 
	\end{eqnarray}
	Taking the derivative of ${\cal F}(s)$ as $s\rightarrow 0$, we have
	\begin{eqnarray}\label{derivative}
	\frac{d{\cal F}}{ds}(s=0)=\lim_{p\rightarrow\infty}\frac{1}{p}\sum_{j=1}^p\min_{\mu_i\in\bR}[g(\mu_j,\beta_{0,j})-\mu_j\hat{\theta}_j]=-\lim_{p\rightarrow\infty}\frac{1}{p}\sum_{j=1}^p\tilde{g}(\hat{\theta}_j,\beta_{0,j}),
	\end{eqnarray}
	where $\hat{\btheta}$ is the solution of (\ref{lse}).
	
	Now we compute $E\log Z_p(\beta,s)$ using the replica method. Instead of performing integration over the log-term directly, we use the following identity
	\begin{eqnarray}
	E\log Z_p(\beta,s)=\left.\frac{\partial E\{Z_p(\beta,s)^d\}}{\partial d}\right|_{d=0}.
	\end{eqnarray}
	We first compute the expectation for integer $d$, and then extrapolate it as $d\rightarrow 0$. 
	
	Define the measure $\nu(d\btheta)$ over $\btheta\in\bR^p$ as follows
	\begin{eqnarray}
	\nu(d\btheta)=\int\exp\left\{-\beta\lambda\|\btheta\|_q^q-\beta s\sum_{j=1}^p[g(\mu_j,\theta_{0,j})-\mu_j\theta_j]\right\}d\bmu d\btheta,
	\end{eqnarray}
	where the integration is with respect to $d\bmu$. Then for integer $d$, define $\nu^d(d\btheta)=\nu(d\btheta^1)\cdots\nu(d\btheta^d)$ be a measure over $(\bR^p)^d$, with $\btheta^1,\cdots,\btheta^d\in\bR^p$. With these notations and substituting (1.1), we have
	\begin{eqnarray}\nn
	E\{Z_p(\beta,s)^d\}&=&\int E\exp\left\{-\frac{\beta}{2n}\sum_{a=1}^d\sum_{i=1}^n\left[\bx_i^T(\bbeta_0-\btheta^a)+w_i\right]^2\right\}\nu^d(d\btheta),
	\end{eqnarray} 
	where the expectation is with respect to $\bw$ and $\bx_i$. Using identity
	\begin{eqnarray}
	\int\delta(u_i^a-\bx_i^T(\bbeta_0-\btheta^a))du_i^a=1,
	\end{eqnarray}
	and Fourier transformation
	\begin{eqnarray}
	\delta(u_i^a-\bx_i^T(\bbeta_0-\btheta^a))=\frac{1}{2\pi i}\int\exp\{i\hat{u}_i^a(u_i^a-\bx_i^T(\bbeta_0-\btheta^a))\}d\hat{u}_i^a,
	\end{eqnarray}
	where $\delta(\cdot)$ is the Dirac delta function, we obtain
	\begin{eqnarray}\nn
	E\{Z_p(\beta,s)^d\}&=&\frac{1}{(2\pi i)^{nd}}\left(\frac{\beta}{n}\right)^{\frac{nd}{2}}\int E\exp\left\{-\frac{\beta}{2n}\sum_{a,i}(u_i^a+w_i)^2\right.\\
	&&\left.+i\sqrt{\frac{\beta}{n}}\sum_{a,i}\hat{u}_i^a(u_i^a-\bx_i^T(\bbeta_0-\btheta^a))\right\}\nu^d(d\btheta)\nu^d(d\bu)\nu^d(d\hat{\bu}),
	\end{eqnarray}
	where $\nu^d(d\bu)=\nu(d\bu^1)\cdots\nu(d\bu^d)$ is a measure over $(\bR^n)^d$ with $\bu^1,\cdots,\bu^d\in\bR^n$ and $\nu^d(d\hat{\bu})$ is defined similarly. For fixed $\sum_{a=1}^d(\bbeta_0-\btheta^a)$, the product term $\bx_i^T\sum_{a=1}^d(\bbeta_0-\btheta^a)$ follows a multivariate normal distribution with mean zero and covariance matrix $\sum_{a,b}(\bbeta_0-\btheta^a)^T\bSigma(\bbeta_0-\btheta^b)$. After integration over $\bx_1,\cdots,\bx_n$, we have
	\begin{eqnarray}\nn
	E\{Z_p(\beta,s)^d\}&=&\frac{1}{(2\pi i)^{nd}}\left(\frac{\beta}{n}\right)^{\frac{nd}{2}}\int E\exp\left\{-\frac{\beta}{2n}\sum_{a,i}(u_i^a+w_i)^2+i\sqrt{\frac{\beta}{n}}\sum_{a,i}\hat{u}_i^au_i^a\right.\\\label{ex}
	&&\left.-\frac{\beta}{2n}\sum_{a,b}(\hat{\bu}^a)^T\hat{\bu}^b(\bbeta_0-\btheta^a)^T\bSigma(\bbeta_0-\btheta^b)\right\}\nu^d(d\btheta)\nu^d(d\bu)\nu^d(d\hat{\bu}).
	\end{eqnarray}
	We next follow the similar procedure as in \cite{Montanari} and use the identity
	\begin{eqnarray}
	\exp(-xy)=\frac{1}{2\pi i}\int_{(-i\infty,i\infty)}\int_{(-\infty,\infty)}\exp(-\zeta \hat{\zeta}+\zeta x-\hat{\zeta}y)d\zeta d\hat{\zeta}
	\end{eqnarray}
	to (\ref{ex}) and introduce integration variables $\bQ=(Q_{ab})_{1\le a,b\le d}$ and $\bLambda=(\Lambda_{ab})_{1\le a,b\le d}$. Letting $d\bQ=\prod_{a,b}dQ_{ab}$  and $d\bLambda=\prod_{a,b}d\Lambda_{ab}$, we have
	\begin{eqnarray}\label{saddle}
	E\{Z_p(\beta,s)^d\}&=&\left(\frac{\beta n}{4\pi i}\right)^{d^2}\int\exp\{-p{\cal S}_d(\bQ,\bLambda)\}d\bQ d\bLambda,\\\label{saddle1}
	{\cal S}_d(\bQ,\bLambda)&=&\frac{\beta\delta}{2}\sum_{a,b}\Lambda_{ab}Q_{ab}-\frac{1}{p}\log\xi(\bLambda)-\delta\log\hat{\xi}(\bQ),\\\nn
	\xi(\bLambda)&=&\int\exp\left\{\frac{\beta}{2}\sum_{a,b}\Lambda_{ab}(\bbeta_0-\btheta^a)^T\bSigma(\bbeta_0-\btheta^b)\right\}\nu^d(d\btheta),\\\nn
	\hat{\xi}(\bQ)&=&\frac{1}{(2\pi i)^d}E\exp\left\{\frac{\beta}{2n}w^2\sum_{a,b}(\bI_{d\times d}+(\beta\bQ)^{-1})^{-1}_{ab}-\frac{1}{2}\log\text{det}(\bI_{d\times d}+\beta\bQ)-\frac{\beta}{2n}dw^2\right\},
	\end{eqnarray}
	where in $\hat{\xi}(\bQ)$ we have performed integrations over $d\bu$ and $d\hat{\bu}$ and the remaining expectation is over the noise variable $w$. We next use the saddle point method in (\ref{saddle}) to obtain
	\begin{eqnarray}
	-\lim_{p\rightarrow\infty}\frac{1}{p}\log E\{Z_p(\beta,s)^d\}={\cal S}_d(\bQ^\star,\bLambda^\star),
	\end{eqnarray}
	where $\bQ^\star,\bLambda^\star$ is the saddle-point location. From replica symmetry, $\bQ^\star,\bLambda^\star$ are invariant under permutations of the row/column indices, which is equivalent to
	\begin{eqnarray}
	Q^\star_{ab}=\left\{\begin{array}{cc}q_1&if~a=b,\\q_0&\text{otherwise}\end{array}\right.,&&
	\Lambda^\star_{ab}=\left\{\begin{array}{cc}\beta\zeta_1&if~a=b,\\\beta\zeta_0&\text{otherwise}\end{array}\right.,
	\end{eqnarray}
	The next step is to substitute the above expressions for $\bQ^\star$ and $\bLambda^\star$ in (\ref{saddle1}) and then taking the limit $d\rightarrow 0$. Hence
	\begin{eqnarray}\nn
	\lim_{d\rightarrow 0}\frac{\beta\delta}{2d}\sum_{a,b}\Lambda^\star_{ab}Q^\star_{ab}&=&\frac{\beta^2\delta}{2}(\zeta_1q_1-\zeta_0q_0)\\\nn
	\lim_{d\rightarrow 0}\frac{-\delta\log\hat{\xi}(\bQ^\star)}{d}&=&\frac{\delta}{2}\log(1+\beta(q_1-q_0))+\frac{\delta}{2}\frac{\beta(q_0+\sigma_w^2)}{1+\beta(q_1-q_0)}\\\nn
	\lim_{d\rightarrow 0}\frac{-\log\xi(\bLambda^\star)}{d}
	&=&E\left\{\log\left[\int\exp\left\{\frac{\beta^2}{2}(\zeta_1-\zeta_0)\|\btheta-\bbeta_0\|_{\bSigma}^2\right.\right.\right.\\\nn
	&&\left.\left.\left.+\beta\sqrt{\zeta_0}\bz^T\bSigma^{1/2}(\btheta-\bbeta_0)\right\}\nu(d\btheta)\right]\right\},
	\end{eqnarray}
	where expectation is with respect to $\bz\sim N(0,\bI_{p\times p})$. In order for the exponent in above equation to be extensive in $p$, we introduce $q_1-q_0=q/\beta$ and $\zeta_1-\zeta_0=-\zeta/\beta$. We find in the asymptotic limit $\beta\rightarrow\infty$ the free energy becomes
	\begin{eqnarray}\nn
	{\cal F}(s)&=&-\lim_{\beta\rightarrow\infty}\lim_{p\rightarrow\infty}\frac{1}{p\beta}E\log Z_p(\beta,s)=\lim_{\beta\rightarrow\infty}\lim_{d\rightarrow 0}\frac{1}{d\beta}{\cal S}_d(\bQ^\star,\bLambda^\star)\\\nn
	&=&\frac{\delta}{2}(\zeta_0q-\zeta q_0)+\frac{\delta}{2}\frac{q_0+\sigma_w^2}{1+q}\\
	&&+\lim_{p\rightarrow\infty}\frac{1}{p}E\min_{\btheta\in \bR^p}\left\{\frac{\zeta}{2}\|\btheta-\bbeta_0\|_{\bSigma}^2-\sqrt{\zeta_0}\bz^T\bSigma^{1/2}(\btheta-\bbeta_0)+D(\btheta;s)\right\},
	\end{eqnarray}
	where
	\begin{eqnarray}
	D(\btheta;s)&=&\min_{\bmu\in\bR^p}\left\{\lambda \|\btheta\|_q^q-s\btheta^{T}\bmu+s\sum_{j=1}^pg(\mu_j,\btheta_{0,j})\right\}.
	\end{eqnarray}
	After eliminating $q$ and $q_0$ using their saddle-point equations and renaming $\zeta_0=\zeta^2\tau^2$, we have
	\begin{eqnarray}\nn
	{\cal F}(s)&=&-\frac{1}{2}(1-\delta)\zeta\tau^2-\frac{\delta}{2}\zeta^2\tau^2+\frac{\delta}{2}\sigma_0^2\zeta\\\label{final}
	&&+\lim_{p\rightarrow\infty}\frac{1}{p}E\min_{\btheta\in \bR^p}\left\{\frac{\zeta}{2}\|\btheta-\bbeta_0-\tau\bSigma^{-1/2}\bz\|_{\bSigma}^2+D(\btheta,s)\right\}.
	\end{eqnarray}
	In the case $s=0$, solving the saddle-point equations for $\zeta$ and $\tau^2$, we finally get
	\begin{eqnarray}\nn
	\tau^2&=&\sigma_w^2+\frac{1}{\delta}\lim_{p\rightarrow\infty}\frac{1}{p\delta}E\left(\|\veta_q(\bbeta_0+\tau\bSigma^{-1/2}\bz,\lambda/\zeta)-\bbeta_0\|_{\bSigma}^2\right),\\
	\zeta&=&1-\frac{1}{\delta}E\{\langle\veta_q^\prime(\bbeta_0+\tau\bSigma^{-1/2}\bz,\lambda/\zeta)\rangle\},
	\end{eqnarray}
	where $\veta_q(\bv,\lambda)$ is defined in (2.2). Then taking the derivative of (\ref{final}) as $s\rightarrow 0$ and minimizing over $\bmu$, we get
	\begin{eqnarray}\label{limit}
	\frac{d{\cal F}}{ds}(s=0)=-\lim_{p\rightarrow\infty}\frac{1}{p}\sum_{j=1}^pE\tilde{g}(\tilde{\theta}_j,\beta_{0,j}),
	\end{eqnarray}
	where $\tilde{\btheta}=\eta_q(\bbeta_0+\tau\bSigma^{-1/2}\bz,\lambda/\zeta)$. Comparing (\ref{limit}) with (\ref{derivative}) for a complete set of functions $\tilde{g}$ and using the standard weak convergence arguments, we prove the claim that the distributional limit of LQLS estimator does indeed hold. 
\end{proof}
 
\subsection{Proof of Proposition \ref{prop1}}\label{prop1p}
\begin{proof}
	Define $\lambda=\alpha_0\tau^{2-q}$. This will enable us to employ the scale invariance properties of the general thresholding operation (\ref{veta}) more efficiently. We can write (\ref{psi}) as
	\begin{eqnarray}\label{psiq}
	\psi_{q}(\tau^2,\alpha_0\tau^{2-q})&=&\frac{1}{p\delta}E\|\veta_q(\bbeta_0+\tau\bSigma^{-1/2}\bz,\alpha_0\tau^{2-q})-\bbeta_0\|^2_{\Sigma}
	\end{eqnarray} 
	Note that $\tau^2=0$ is actually a fixed point of $\psi_{q}(\tau^2,\alpha_0\tau^{2-q})$. Furthermore, it is straightforward to see that 0 is a stable fixed point if and only if
	\begin{eqnarray}\nn
	\frac{d\psi_{q}(\tau^2,\alpha_0\tau^{2-q})}{d\tau^2}\left|_{\tau^2=0}=\lim_{\tau^2\rightarrow 0}\frac{\psi_{q}(\tau^2,\alpha_0\tau^{2-q})}{\tau^2}\textless 1\right.
	\end{eqnarray} 
	Denote $U=\bbeta_0/\tau$. Then (\ref{psiq}) can be written as
	\begin{eqnarray}\nn
	\psi_{q}(\tau^2,\alpha_0\tau^{2-q})&=&\frac{\tau^2}{\delta}R_q(\tau^2,\alpha_0),
	\end{eqnarray} 
	where
	\begin{eqnarray}\label{rfun}
	R_q(\tau^2,\alpha_0)&=&\frac{1}{p}E\|\hat{\bbeta}-U\|^2_{\Sigma},
	\end{eqnarray} 
	where
	\begin{eqnarray}\label{betahat}
	\hat{\bbeta}&=&\text{argmin}_{\bbeta}\left\{\frac{1}{2}\|\bbeta-\bU-\Sigma^{-1/2}\bz\|_\Sigma^2+\alpha_0\|\bbeta\|_q^q\right\}.
	\end{eqnarray} 
	For fixed $U$, denote $S=\{j|U_{j}\ne 0\}$. As $\tau\rightarrow 0$, $\hat{\bbeta}_S=U_s+o_P(\mu/\tau)$, we obtain
	\begin{eqnarray}\nn
	\bSigma_{SS}(\hat{\bbeta}_S-\bU_S)+\bSigma_{S\bar{S}}\hat{\bbeta}_{\bar{S}}-(\bSigma^{1/2}\bz)_S+q\alpha_0|\hat{\bbeta}_S|^{q-1}\text{sign}(\hat{\bbeta}_S)=0,
	\end{eqnarray} 
	which can be written as
	\begin{eqnarray}\label{b2}
	\hat{\bbeta}_S-\bU_S=(\bSigma^{-1/2}\bz)_S-\bSigma_{SS}^{-1}\bSigma_{S\bar{S}}\hat{\bbeta}_{\bar{S}}-q\alpha_0\bSigma_{SS}^{-1}|\hat{\bbeta}_S|^{q-1}\text{sign}(\hat{\bbeta}_S).
	\end{eqnarray} 
	Substituting into (\ref{rfun}), we have
	\begin{eqnarray}\nn
	R_{q}(\tau^2,\alpha_0)&=&\frac{1}{p}E\{(\hat{\bbeta}_S-\bU_S)^T\bSigma_{SS}(\hat{\bbeta}_S-\bU_S)+2(\hat{\bbeta}_S-\bU_S)^T\bSigma_{S\bar{S}}\hat{\bbeta}_{\bar{S}}+\hat{\bbeta}_{\bar{S}}^T\bSigma_{\bar{S}\bar{S}}\hat{\bbeta}_{\bar{S}}\}\\\nn
	&=&\frac{1}{p}E[\{(\bSigma^{1/2}\bz)_S-q\alpha_0|\hat{\bbeta}_S|^{q-1}\text{sign}(\hat{\bbeta}_S)+\bSigma_{S\bar{S}}\hat{\bbeta}_{\bar{S}}\}^T\bSigma_{SS}^{-1}\\\nn
	&&\{(\bSigma^{1/2}\bz)_S-q\alpha_0|\hat{\bbeta}_S|^{q-1}\text{sign}(\hat{\bbeta}_S)-\bSigma_{S\bar{S}}\hat{\bbeta}_{\bar{S}}\}+\hat{\bbeta}_{\bar{S}}^T\bSigma_{\bar{S}\bar{S}}\hat{\bbeta}_{\bar{S}}]\\\nn
	&=&\frac{1}{p}E[\{(\bSigma^{1/2}\bz)_S-q\alpha_0|\hat{\bbeta}_S|^{q-1}\text{sign}(\hat{\bbeta}_S)\}^T\bSigma_{SS}^{-1}\{(\bSigma^{1/2}\bz)_S-q\alpha_0|\hat{\bbeta}_S|^{q-1}\text{sign}(\hat{\bbeta}_S)\}\\\label{rqd}
	&&+\hat{\bbeta}_{\bar{S}}^T(\bSigma_{\bar{S}\bar{S}}-\bSigma_{\bar{S}S}\bSigma_{SS}^{-1}\bSigma_{S\bar{S}})\hat{\bbeta}_{\bar{S}}]
	\end{eqnarray} 
	As $\tau\rightarrow 0$, $\hat{\bbeta}_S\rightarrow U_S$. Thus, $|\hat{\bbeta}_S|^{q-1}=O_p(\tau^{1-q})$, we have
	\begin{eqnarray}\nn
	R_{q}(\tau^2,\alpha_0)&=&\frac{1}{p}E[(\bSigma^{1/2}\bz)_S^T\bSigma_{SS}^{-1}(\bSigma^{1/2}\bz)_S+\hat{\bbeta}_{\bar{S}}^T(\bSigma_{\bar{S}\bar{S}}-\bSigma_{\bar{S}S}\bSigma_{SS}^{-1}\bSigma_{S\bar{S}})\hat{\bbeta}_{\bar{S}}]+O(\tau^{1-q}).
	\end{eqnarray} 
	The second term $\hat{\bbeta}_{\bar{S}}^T(\bSigma_{\bar{S}\bar{S}}-\bSigma_{\bar{S}S}\bSigma_{SS}^{-1}\bSigma_{S\bar{S}})\hat{\bbeta}_{\bar{S}}\ge 0$. Combining (\ref{betahat}) and (\ref{b2}), we obtain the equation for $\hat{\bbeta}_{\bar{S}}$ as
	\begin{eqnarray}\nn
	\hat{\bbeta}_{\bar{S}}=\text{argmin}_{\bbeta}\left\{|\hat{\bbeta}-\overline{\bSigma}_{\bar{S}\bar{S}}^{-1}\{(\bSigma^{1/2}\bz)_S-\bSigma_{\bar{S}S}(\bSigma^{1/2}\bz)_{\bar{S}}\}|_{\overline{\bSigma}_{\bar{S}\bar{S}}}+q\alpha_0|\hat{\bbeta}|^{q-1}\text{sign}(\hat{\bbeta})\right\}
	\end{eqnarray} 
	which goes to 0 if $\alpha_0$ is large enough. Here $\overline{\bSigma}_{\bar{S}\bar{S}}= \bSigma_{\bar{S}\bar{S}}-\bSigma_{\bar{S}S}\bSigma^{-1}_{SS}\bSigma_{S\bar{S}}$. Therefore denote $\alpha_{0\star}$ the optimal $\alpha_0$, we obtain
	\begin{eqnarray}\nn
	R_{q}(\tau^2,\alpha_{0\star})&\xrightarrow{\tau\rightarrow 0}&\frac{1}{p}E[(\bSigma^{1/2}\bz)_S^T\bSigma_{SS}^{-1}(\bSigma^{1/2}\bz)_S]=\epsilon,
	\end{eqnarray} 
	and 
	\begin{eqnarray}\nn
	\frac{d\psi_{q}(\tau^2,\alpha_{0\star}\tau^{2-q})}{d\tau^2}\left|_{\tau^2=0}=\frac{\epsilon}{\delta}\right..
	\end{eqnarray} 
\end{proof}

\subsection{Proof of Proposition \ref{prop3}}
\begin{proof}	
	Define $\lambda=\alpha_0\tau^2$, then from (\ref{veta}), we have for $q=0$
	\begin{eqnarray}\nn
	\bar{\bbeta}&=&\hat{\bbeta}/\tau=\veta_0(\bbeta_0/\tau+\bSigma^{-1/2}\bz,\alpha_0)\\\label{sveta}
	&=&\text{argmin}_{\bbeta\in\bR^p}\left\{\frac{1}{2}\|\bbeta-\bbeta_0/\tau-\bSigma^{-1/2}\bz\|^2_{\bSigma}+\alpha_0\|\bbeta\|_0\right\}.
	\end{eqnarray} 
	For block diagonal covariance $\bSigma$ with block $\left(\begin{array}{cc}1&\rho\\\rho&1\end{array}\right)$, equation (\ref{sveta}) can be decomposed into 2-dimensional blocks each of which can be solved as    
	\begin{eqnarray}\nn
	(\bar{\beta}_1,\bar{\beta}_2)=(\hat{\beta}_1/\tau,\hat{\beta}_2/\tau)=\text{argmin}_{\beta_1,\beta_2}{\cal L},
	\end{eqnarray} 
	where
	\begin{eqnarray}\nn
	{\cal L}&=&\frac{1}{2}\{(\beta_1-y_1)^2+(\beta_2-y_2)^2+2\rho(\beta_1-y_1)(\beta_2-y_2)\}+\alpha_0(\|\beta_1\|_0+\|\beta_2\|_0).
	\end{eqnarray} 
	For problem (\ref{sveta}), $y_1=U_1+\xi_1$, $y_2=U_2+\xi_2$, where $U_1=\beta_{0,1}/\tau$, $U_2=\beta_{0,2}/\tau$, $\beta_{0,1},\beta_{0,2}\overset{i.i.d.}{\sim} G$, and 
	\begin{eqnarray}\nn
	\xi_1&=&\frac{1}{2}\left(\sqrt{\frac{1}{1+\rho}}+\sqrt{\frac{1}{1-\rho}}\right)z_1+\frac{1}{2}\left(\sqrt{\frac{1}{1+\rho}}-\sqrt{\frac{1}{1-\rho}}\right)z_2,\\\nn 
	\xi_2&=&\frac{1}{2}\left(\sqrt{\frac{1}{1+\rho}}-\sqrt{\frac{1}{1-\rho}}\right)z_1+\frac{1}{2}\left(\sqrt{\frac{1}{1+\rho}}+\sqrt{\frac{1}{1-\rho}}\right)z_2,
	\end{eqnarray} 
	with $z_1,z_2\overset{i.i.d.}{\sim}N(0,1)$. 
	
	As shown in Figure \ref{figure1}, the two dimensional space can be divided into nine regions. In each region, the estimator is chosen to be the one that gives the lowest loss ${\cal L}$. The solution can be summarized as
	\begin{eqnarray}\label{ninearea}
	\left\{
	\begin{array}{ccc}
	\hat{\beta}_1\textgreater 0~\&~\hat{\beta}_2\textgreater 0&if&\by\in I_1\\
	\hat{\beta}_1\textgreater 0~\&~\hat{\beta}_2=0&if&\by\in I_5\\
	\hat{\beta}_1\textgreater 0~\&~\hat{\beta}_2\textless 0&if&\by\in I_4\\
	\hat{\beta}_1= 0~\&~\hat{\beta}_2\textgreater 0&if&\by\in I_7\\
	\hat{\beta}_1=0~\&~\hat{\beta}_2\textless 0&if&\by\in I_8\\
	\hat{\beta}_1\textless 0~\&~\hat{\beta}_2\textgreater 0&if&\by\in I_2\\
	\hat{\beta}_1\textless 0~\&~\hat{\beta}_2= 0&if&\by\in I_6\\
	\hat{\beta}_1\textless 0~\&~\hat{\beta}_2\textless 0&if&\by\in I_3\\
	\hat{\beta}_1=0~\&~\hat{\beta}_2=0&if&\by\in I_9 
	\end{array}
	\right.
	\end{eqnarray} 
	where
	\begin{eqnarray}\nn
	\left\{
	\begin{array}{ccl}
	I_1&=&I\left(y_1\textgreater\sqrt{\frac{2\alpha_0}{1-\rho^2}}~\&~y_2\textgreater\sqrt{\frac{2\alpha_0}{1-\rho^2}}\right)\\
	I_5&=&I\left(|y_2|\textless\sqrt{\frac{2\alpha_0}{1-\rho^2}}~\&~|y_1|\textgreater| y_2|~\&~y_1+\rho y_2\textgreater\sqrt{2\alpha_0}\right)\\
	I_4&=&I\left(y_1\textgreater\sqrt{\frac{2\alpha_0}{1-\rho^2}}~\&~y_2\textless-\sqrt{\frac{2\alpha_0}{1-\rho^2}}~\&~y_1^2+y_2^2+2\rho y_1y_2\textgreater\sqrt{4\alpha_0}\right)\\
	I_7&=&I\left(|y_1|\textless\sqrt{\frac{2\alpha_0}{1-\rho^2}}~\&~|y_1|\textless|y_2|~\&~y_2+\rho y_1\textgreater\sqrt{2\alpha_0}\right)\\
	I_8&=&I\left(|y_1|\textless\sqrt{\frac{2\alpha_0}{1-\rho^2}}~\&~|y_1|\textless|y_2|~\&~y_2+\rho y_1\textless -\sqrt{2\alpha_0}\right)\\
	I_2&=&I\left(y_1\textless-\sqrt{\frac{2\alpha_0}{1-\rho^2}}~\&~y_1^2+y_2^2+2\rho y_1y_2\textgreater\sqrt{4\alpha_0}~\&~y_2\textgreater\sqrt{\frac{2\alpha_0}{1-\rho^2}}\right)\\
	I_6&=&I\left(|y_2|\textless\sqrt{\frac{2\alpha_0}{1-\rho^2}}~\&~|y_1|\textgreater|y_2|~\&~y_1+\rho y_2\textless -\sqrt{2\alpha_0}\right)\\
	I_3&=&I\left(y_1\textless -\sqrt{\frac{2\alpha_0}{1-\rho^2}}~\&~y_2\textless -\sqrt{\frac{2\alpha_0}{1-\rho^2}}\right)\\
	I_9&=&(\cup_{i=1}^8I_i)^c
	\end{array}
	\right.
	\end{eqnarray} 
	We consider three different scenarios. The first scenario includes regions 1-4 where both $\hat{\beta}_1$ and $\hat{\beta}_2$ are nonzero and the loss is ${\cal L}=2\alpha_0$. Its contribution to the prediction risk is  
	\begin{eqnarray}\nn
	R_1(\tau^2,\alpha_0)&=&E\{(\hat{\beta}_1-\beta_{0,1})^2+(\hat{\beta}_2-\beta_{0,2})^2+2\rho(\hat{\beta}_1-\beta_{0,1})(\hat{\beta}_2-\beta_{0,2})\}(I_1+I_2+I_3+I_4)/2\\\label{r01}
	&=&\tau^2E\{(\xi_1^2+\xi_2^2+2\rho\xi_1\xi_2)\}(I_1+I_2+I_3+I_4)/2.
	\end{eqnarray} 
	The second scenario includes regions 5-8 where $\hat{\beta}_1\ne 0,\hat{\beta}_2=0$ or $\hat{\beta}_1=0,\hat{\beta}_2\ne 0$ and the loss is ${\cal L}=\frac{1}{2}(1-\rho^2)y_2^2+\alpha_0$ or ${\cal L}=\frac{1}{2}(1-\rho^2)y_1^2+\alpha_0$. Its 
	contribution to the prediction risk is  
	\begin{eqnarray}\nn
	R_2(\tau^2,\alpha_0)&=&\tau^2\left[E\{(\xi_1+\rho\xi_2)^2+(1-\rho^2)U_2^2\}(I_5+I_6)\right.\\\label{r02}
	&&\left.+E\{(\xi_2+\rho\xi_1)^2+(1-\rho^2)U_1^2\}(I_7+I_8)\right]/2.
	\end{eqnarray} 
	The third scenario includes regions 9 where $\hat{\beta}_1=\hat{\beta}_2=0$ and the loss is ${\cal L}=(y_1^2+y_2^2+2\rho y_1y_2)/2$. Its contribution to the prediction risk is  
	\begin{eqnarray}\label{r03}
	R_3(\tau^2,\alpha_0)&=&\tau^2E\{(U_1^2+U_2^2+2\rho U_1U_2)\}I_9/2.
	\end{eqnarray} 
	
	\begin{figure}[hbtp]
		\begin{center}
			\epsfig{file=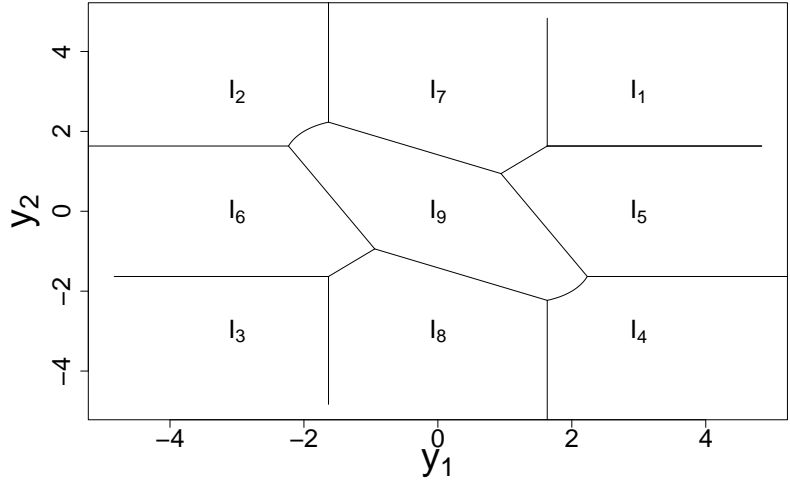,width=8.1cm,angle=-90}
		\end{center}
		\caption{Illustration of the solution (\ref{ninearea}) for equation (\ref{sveta}) in two dimensional space. Here $\rho=0.5$ and $\alpha_0=1$.}
		\label{figure1}
	\end{figure}
	
	To find the optimal $\alpha_0$ for minimizing $R(\tau^2,\alpha_0)=R_1(\tau^2,\alpha_0)+R_2(\tau^2,\alpha_0)+R_3(\tau^2,\alpha_0)$, we need to take derivative over $\alpha_0$. Since the explicit function $R(\tau^2,\alpha_0)$ in each region does not depend on $\alpha_0$, changing $\alpha_0$ only causes the change of the boundaries among different regions. According to Stokes's theorem, as in Theorem 1 of \cite{baddeley_1977}, we conclude that the dominant contributions to $\lim_{\tau\rightarrow 0}\frac{\partial R(\tau^2,\alpha_0)}{\partial\alpha_0}$ come from the boundaries between scenarios 1 and 2 as well as the boundaries between scenarios 2 and 3. For example, the contribution from the boundary between regions 1 and 5 is
	\begin{eqnarray}\nn
	&&\tau^2\sqrt{\frac{1-\rho^2}{2\alpha_0}}E\left\{\left(\sqrt{\frac{2\alpha_0}{1-\rho^2}}-U_2\right)^2-U_2^2\right\}I\left(U_1+\xi_1\textgreater\sqrt{\frac{2\alpha_0}{1-\rho^2}}\right)\delta\left(U_2+\xi_2=\sqrt{\frac{2\alpha_0}{1-\rho^2}}\right)\\\nn
	&\sim&\epsilon_+^2E\left\{\left(\sqrt{\frac{2\alpha_0}{1-\rho^2}}-U_2\right)^2-U_2^2\right\}\phi(\sqrt{2\alpha_0}-\sqrt{1-\rho^2}U_2)+\epsilon_+(1-\epsilon)\frac{2\alpha_0}{1-\rho^2}\phi(\sqrt{2\alpha_0}),
	\end{eqnarray}
	where $\epsilon_+=P(\beta_0\textgreater 0)$. The contribution from the boundary between regions 5 and 9 is 
	\begin{eqnarray}\nn
	&&\tau^2\sqrt{\frac{1}{2\alpha_0}}E\{(\sqrt{2\alpha_0}-U_1+\rho U_2)^2-(U_1+\rho U_2)^2\}I\left(\frac{\sqrt{2\alpha_0}}{1+\rho}\textgreater U_2+\xi_2\textgreater-\sqrt{\frac{2\alpha_0}{1-\rho^2}}\right)\\\nn
	&&\delta\left(U_1+\rho U_2+\xi_1+\rho\xi_2=\sqrt{2\alpha_0}\right)\\\nn
	&\sim&(1-\epsilon)^2E\alpha_0\phi(\sqrt{2\alpha_0}).
	\end{eqnarray} 
	The contributions from other boundaries can be derived similarly. The total contributions include three dominant terms: $\phi(\sqrt{2\alpha_0})$, $\phi(\sqrt{2\alpha_0}-\sqrt{1-\rho^2}U_1)$, and $\phi(\sqrt{2\alpha_0}-\sqrt{1-\rho^2}U_2)$. Therefore, we need to have $\lim_{\tau\rightarrow 0}\sqrt{2\alpha_{0\star}}\tau=\mu\sqrt{1-\rho^2}/2$ to ensure the existence of the finite solutions for $\frac{\partial R(\tau^2,\alpha_0)}{\partial\alpha_0}=0$. Substituting into (\ref{r01}), (\ref{r02}), and (\ref{r03}), we obtain	
	\begin{eqnarray}\nn
	R_1(\tau^2,\alpha_{0\star})&=&\frac{\tau^2}{2}E\{(\xi_1^2+\xi_2^2+2\rho\xi_1\xi_2)\}I\left(|U_1+\xi_1|\textgreater\sqrt{\frac{2\alpha_{0\star}}{1-\rho^2}}\right)I\left(|U_2+\xi_2|\textgreater\sqrt{\frac{2\alpha_{0\star}}{1-\rho^2}}\right)\\\nn
	&=&\tau^2[\epsilon^2(1+o(\sqrt{2\alpha_{0\star}}\phi(\sqrt{2\alpha_{0\star}})))+\epsilon(1-\epsilon)o(\sqrt{2\alpha_{0\star}}\phi(\sqrt{2\alpha_{0\star}}))+(1-\epsilon)^2o(\phi(\sqrt{2\alpha_{0\star}})^2)]\\\nn
	&=&\tau^2[\epsilon^2+o(\sqrt{2\alpha_{0\star}}\phi(\sqrt{2\alpha_{0\star}}))],
	\end{eqnarray} 
	\begin{eqnarray}\nn
	R_2(\tau^2,\alpha_{0\star})&=&\frac{\tau^2}{2}E\{(\xi_1+\rho\xi_2)^2+(1-\rho^2)U_2^2\}I\left(|y_1+\rho y_2|\textgreater\sqrt{2\alpha_{0\star}}\right)I\left(|U_2+\xi_2|\textless\sqrt{\frac{2\alpha_{0\star}}{1-\rho^2}}\right)\\\nn
	&&+\frac{\tau^2}{2}E\{(\xi_2+\rho\xi_1)^2+(1-\rho^2)U_1^2\}I\left(|y_2+\rho y_1|\textgreater\sqrt{2\alpha_{0\star}}\right)I\left(|U_1+\xi_1|\textless\sqrt{\frac{2\alpha_{0\star}}{1-\rho^2}}\right)\\\nn
	&=&\tau^2[\epsilon^2o(\sqrt{2\alpha_{0\star}}\phi(\sqrt{2\alpha_{0\star}}))+\epsilon(1-\epsilon)(1+o(\sqrt{2\alpha_{0\star}}\phi(\sqrt{2\alpha_{0\star}})))+(1-\epsilon)^2o(\sqrt{2\alpha_{0\star}}\phi(\sqrt{2\alpha_{0\star}}))]\\\nn
	&=&\tau^2[\epsilon(1-\epsilon)+o(\sqrt{2\alpha_{0\star}}\phi(\sqrt{2\alpha_{0\star}}))],
	\end{eqnarray} 
	and
	\begin{eqnarray}\nn
	R_3(\tau^2,\alpha_{0\star})&=&\frac{\tau^2}{2}E(U_1^2+U_2^2+2\rho U_1U_2)\left(|y_1+\rho y_2|\textgreater\sqrt{2\alpha_{0\star}}\right)\left(|y_2+\rho y_1|\textgreater\sqrt{2\alpha_{0\star}}\right)\\\nn
	&\le&\frac{\tau^2}{2}E(U_1^2+U_2^2+2\rho U_1U_2)I\left(|U_1+\xi_1|\textless\frac{(\rho+\sqrt{1-\rho^2})\sqrt{2\alpha}}{\sqrt{1-\rho^2}}\right)\\\nn
	&&I\left(|U_2+\xi_2|\textless\frac{(\rho+\sqrt{1-\rho^2})\sqrt{2\alpha_{0\star}}}{\sqrt{1-\rho^2}}\right)\\\nn
	&=&\tau^2[\epsilon^2o(\phi((2-\rho-\sqrt{1-\rho^2})\sqrt{2\alpha_{0\star}})^2)+\epsilon(1-\epsilon)o(\sqrt{2\alpha_{0\star}}\phi((2-\rho-\sqrt{1-\rho^2})\sqrt{2\alpha_{0\star}}))]\\\nn
	&=&\tau^2o(\sqrt{2\alpha_{0\star}}\phi((2-\rho-\sqrt{1-\rho^2})\sqrt{2\alpha_{0\star}})).
	\end{eqnarray} 
	Thus, as $\tau\rightarrow 0$, the total MSE becomes $R(\tau^2,\alpha_{0\star})=\tau^2(\epsilon+o(\phi(\tilde{\mu}\tau^{-1}))$, where $\tilde{\mu}=(2-\rho-\sqrt{1-\rho^2})\sqrt{1-\rho^2}\mu/2$.
	The fixed point equation for small $\sigma_w^2$ can be written as
	\begin{eqnarray}\nn
	\tau_l^2&=&\sigma_w^2+\frac{R(\tau_l^2,\alpha_{0\star})}{\delta}\\\nn
	&=&\sigma_w^2+\frac{\tau^2_l}{\delta}(\epsilon+o(\phi(\tilde{\mu}\tau_l^{-1}))).
	\end{eqnarray} 
	Therefore,
	\begin{eqnarray}\nn
	\tau_l^2&=&\frac{\delta\sigma_w^2}{\delta-\epsilon}+o\left(\phi\left(\tilde{\mu}\sqrt{\frac{\delta-\epsilon}{\delta}}\sigma_w^{-1}\right)\right).
	\end{eqnarray}
	
\end{proof}

\subsection{Proof of Proposition \ref{prop4}}
\begin{proof}
	For $\bSigma=\left(\begin{array}{cc}1&\rho\\\rho&1\end{array}\right)$, denote $\rho_1=\frac{1}{2}(\sqrt{1+\rho}+\sqrt{1-\rho})$ and $\rho_2=\frac{1}{2}(\sqrt{1+\rho}-\sqrt{1-\rho})$. Let $\lambda=\alpha_0\tau^{2-q}$, then the solution can be decomposed into 2-dimensional block of $\bbeta\in\bR^2$ as
	\begin{eqnarray}\nn
	\hat{\beta}_1,\hat{\beta}_2=\text{argmin}_{\beta_1,\beta_2}\left[{\cal L}\right],
	\end{eqnarray} 
	where
	\begin{eqnarray}\nn
	{\cal L}&=&\frac{1}{2}\|\bbeta-\bbeta_0-\tau\Sigma^{-1/2}\bz\|_\Sigma^2+\alpha_0\tau^{2-q}\|\bbeta\|_q^q\\\nn
	&=&\frac{1}{2}\{(\beta_1-\beta_{0,1})^2+(\beta_2-\beta_{0,2})^2+2\rho(\beta_1-\beta_{0,1})(\beta_2-\beta_{0,2})\}\\\label{lfun1}
	&&-\tau\xi_1(\beta_1-\beta_{0,1})-\tau\xi_2(\beta_2-\beta_{0,2})+\frac{1}{2}\tau^2\|\bz\|^2+\alpha_0\tau^{2-q}(|\beta_1|^q+|\beta_2|^q),
	\end{eqnarray} 
	with $\xi_1=\rho_1z_1+\rho_2z_2$, $\xi_2=\rho_2z_1+\rho_1z_2$. If both $\hat{\beta}_1$ and $\hat{\beta}_2$ are nonzero, they should satisfy 
	\begin{eqnarray}\label{2deq}
	\left\{\begin{array}{ccc}
	\hat{\beta}_1-\beta_{0,1}+\rho(\hat{\beta}_2-\beta_{0,2})-\tau\xi_1+\alpha_0\tau^{2-q}q|\hat{\beta}_1|^{q-1}\text{sign}(\hat{\beta}_1)&=&0,\\
	\rho(\hat{\beta}_1-\beta_{0,1})+\hat{\beta}_2-\beta_{0,2}-\tau\xi_2+\alpha_0\tau^{2-q}q|\hat{\beta}_2|^{q-1}\text{sign}(\hat{\beta}_2)&=&0.
	\end{array}\right.
	\end{eqnarray} 
	The $\psi$ function defined in (2.3) can now be written as
	\begin{eqnarray}\nn
	\psi_q(\tau^2,\alpha_0\tau^{2-q})&=&\sigma_w^2+\frac{1}{2\delta}E\|\hat{\bbeta}-\bbeta_0\|_\Sigma\\\label{2dpsi}
	&=&\sigma_w^2+\frac{1}{2\delta}E\{(\hat{\beta}_1-\beta_{0,1})^2+(\hat{\beta}_2-\beta_{0,2})^2+2\rho(\hat{\beta}_1-\beta_{0,1})(\hat{\beta}_2-\beta_{0,2})\}.
	\end{eqnarray} 
	We consider three different situations. In the first case: $\beta_{0,1}\ne 0$ and $\beta_{0,2}\ne 0$. Then as $\tau\rightarrow 0$, both $\hat{\beta}_1$ and $\hat{\beta}_2$ are nonzero and (\ref{2deq}) can be written as
	\begin{eqnarray}\nn
	\hat{\bbeta}&=&\bbeta_0+\tau\Sigma^{-1/2}\bz-\alpha_0\tau^{2-q}q\bSigma^{-1}|\bbeta_0|^{q-1}\text{sign}(\bbeta_0)+o_p(\tau^{2-q}).
	\end{eqnarray} 
	Substituting into (\ref{2dpsi}), its contribution to the prediction risk is
	\begin{eqnarray}\nn
	R_1(\tau^2,\alpha_0\tau^{2-q})&=&\epsilon^2\tau^2E\left[1+\frac{1}{2}\alpha_0^2\tau^{2-2q}q^2\{(\bbeta_0)^{q-1}sign(\bbeta_0)\}^T\bSigma^{-1}(\bbeta_0)^{q-1}sign(\bbeta_0)\right]\\\label{r10}
	&&+o(\tau^{4-2q}).
	\end{eqnarray} 
	Taking derivative over $\alpha_0$, we obtain
	\begin{eqnarray}\nn
	\frac{\partial R_1(\tau^2,\alpha_0\tau^{2-q})}{\tau^2\partial\alpha_0}&=&\epsilon^2E\left[\alpha_0\tau^{2-2q}q^2\{(\bbeta_0)^{q-1}sign(\bbeta_0)\}^T\bSigma^{-1}|\bbeta_0|^{q-1}sign(\bbeta_0)\right]+o(\tau^{2-2q})\\\label{r1}
	&=&2\epsilon^2\alpha_0\tau^{2-2q}q^2D_0+o(\tau^{2-2q}).
	\end{eqnarray} 
	where $D_0=E\{|\beta_0|^{2q-2}\}-\rho(E\{|\beta_0|^{q-1}sign(\beta_0)\})^2/(1-\rho^2)$.
	
	In the second case: $\beta_{0,1}\ne 0$ and $\beta_{0,2}=0$ (or $\beta_{0,2}\ne 0$ and $\beta_{0,1}=0$). The first equation of (\ref{2deq}) becomes
	\begin{eqnarray}\nn
	\hat{\beta}_1-\beta_{0,1}&=&\tau\xi_1-\rho\hat{\beta}_2-\alpha_0\tau^{2-q}q|\beta_{0,1}|^{q-1}\text{sign}(\beta_{0,1})+o_p(\tau^{2-q}).
	\end{eqnarray} 
	Substituting into (\ref{lfun1}), we obtain the solution for $\beta_2$
	\begin{eqnarray}\nn
	\hat{\beta}_2=\text{argmin}_\beta\left\{(1-\rho^2)[\beta-(1-\rho^2)^{-1}\tau(\xi_2-\rho\xi_1)]^2+\alpha_0\tau^{2-q}q|\beta|^{q-1}\text{sign}(\beta)\right\}
	\end{eqnarray}
	which is zero when $|\xi_2-\rho\xi_1|\textless(1-\rho^2)C_q(\frac{\alpha_0}{1-\rho^2})^{1/(2-q)}$ with $C_q=[2(1-q)]^{1/(2-q)}+q[2(1-q)]^{(q-1)/(2-q)}$ and nonzero solved otherwise Denote $I_1=I\left(|\xi_2-\rho\xi_1|\textgreater(1-\rho^2)C_q(\frac{\alpha_0}{1-\rho^2})^{\frac{1}{2-q}}\right)$. Its contribution to the prediction risk can be written as
	\begin{eqnarray}\nn
	R_2(\tau^2,\alpha_0\tau^{2-q})&=&\epsilon(1-\epsilon)E\{(\tau\xi_1-\alpha_0\tau^{2-q}q|\beta_{0,1}|^{q-1}\text{sign}(\beta_{0,1}))^2+(1-\rho^2)\hat{\beta}_2^2\}\\\label{r20}
	&=&\epsilon(1-\epsilon)\tau^2\{1+\alpha_0^2\tau^{2-2q}q^2E|\beta_{0,1}|^{2q-2}+(1-\rho^2)E\tilde{\beta}^2I_1\}+o(\tau^{2-2q}).
	\end{eqnarray} 
	where $\tilde{\beta}$ is the solution of
	\begin{eqnarray}\nn
	(1-\rho^2)\tilde{\beta}-(\xi_2-\rho\xi_1)+\alpha_0 q|\tilde{\beta}|^{q-1}\text{sign}(\tilde{\beta})=0.
	\end{eqnarray}
	Taking derivative over $\alpha_0$, we obtain
	\begin{eqnarray}\nn
	\frac{\partial R_2(\tau^2,\alpha_0\tau^{2-q})}{\tau^2\partial\alpha_0}&=&2\epsilon(1-\epsilon)\left\{\alpha_0\tau^{2-2q}q^2C_0-(1-\rho^2)D_q^2\left(\frac{\alpha_0}{1-\rho^2}\right)^{\frac{2}{2-q}}C_q\left(\frac{1}{1-\rho^2}\right)^{\frac{q}{2(2-q)}}\frac{1}{2-q}\alpha_0^\frac{q-1}{2-q}\right.\\\label{r2}
	&&~~~~~~~~~~~~~~~~~~~~~~~~~~~~\left.\phi\left(C_q\left(\frac{1}{1-\rho^2}\right)^{\frac{q}{2(2-q)}}\alpha_0^{\frac{1}{2-q}}\right)\right\}+o(\tau^{2-q})\\\nn
	&=&2\epsilon(1-\epsilon)\left\{\alpha_0\tau^{2-2q}q^2C_0-\left(\frac{1}{1-\rho^2}\right)^\frac{3q}{2(2-q)}\frac{C_qD_q^2}{2-q}\alpha_0^\frac{q+1}{2-q}\phi\left(C_q\left(\frac{1}{1-\rho^2}\right)^{\frac{q}{2(2-q)}}\alpha_0^{\frac{1}{2-q}}\right)\right\},
	\end{eqnarray} 
	where $C_0=E|\beta_{0,1}|^{2q-2}$ and $D_q=[2(1-q)]^{\frac{1}{2-q}}$.
	
	In the third case, $\beta_{0,1}=\beta_{0,2}=0$, its contribution to MSE can be summarized as 
	\begin{eqnarray}\label{r30}
	R_3(\tau^2,\alpha_0\tau^{2-q})&=&\frac{1}{2}E(\hat{\beta}_1^2+\hat{\beta}_2^2+2\rho\hat{\beta}_1\hat{\beta}_2),
	\end{eqnarray}
	where $\hat{\beta}_1$ and $\hat{\beta}_2$ are the solution that minimizes
	\begin{eqnarray}\label{lfun0}
	{\cal L}&=&\frac{1}{2}\|\bbeta-\tau\Sigma^{-1/2}\bz\|_\Sigma^2+\alpha_0\tau^{2-q}\|\bbeta\|_q^q.
	\end{eqnarray} 
	As shown in Figure \ref{figure2}, we can divide the $2$-dimensional space into nine regions. In region $I_9$, $\hat{\beta}_1=0$ and $\hat{\beta}_2=0$; in regions $I_5$ and $I_6$, $\hat{\beta}_2=0$ and $\hat{\beta}_1$ is solved by $\hat{\beta}_1-\tau\xi_1+\alpha_0\tau^{2-q} q|\hat{\beta}_1|^{q-1}\text{sign}(\hat{\beta}_1)=0$; in regions $I_7$ and $I_8$, $\hat{\beta}_1=0$ and $\hat{\beta}_2$ is solved by $\hat{\beta}_2-\tau\xi_2+\alpha_0\tau^{2-q} q|\hat{\beta}_1|^{q-1}\text{sign}(\hat{\beta}_2)=0$; In regions $I_1$, $I_2$, $I_3$, and $I_4$, $\hat{\beta}_1\ne 0$ and $\hat{\beta}_2\ne 0$ which are solved by 
	\begin{eqnarray}\nn
	\left\{\begin{array}{ccc}
	\hat{\beta}_1+\rho\hat{\beta}_2-\tau\xi_1+\alpha_0\tau^{2-q} q|\hat{\beta}_1|^{q-1}\text{sign}(\hat{\beta}_1)&=&0,\\
	\rho\hat{\beta}_1+\hat{\beta}_2-\tau\xi_2+\alpha_0\tau^{2-q} q|\hat{\beta}_2|^{q-1}\text{sign}(\hat{\beta}_2)&=&0.
	\end{array}\right.
	\end{eqnarray} 
	Then the derivative of $R_3(\tau^2,\alpha_0\tau^{2-q})$ over $\alpha_0$ involves the explicit derivative inside each region and integrals over the boundaries among different regions over $1$-dimensional boundary curve. According to Stokes's theorem, as in Theorem 1 of \cite{baddeley_1977}, we conclude that the main contributions come from the four dominated boundaries which connect region 9 and four other regions 5, 6, 7, and 8 respectively. The contribution can be summarized as
	\begin{eqnarray}\label{r53}
	\frac{\partial R_3(\tau^2,\alpha_0\tau^{2-q})}{\tau^2\partial\alpha_0}&=&-2(1-\epsilon)^2\frac{C_qD_q^2}{2-q}\alpha_0^{(q+1)/(2-q)}\phi(C_q\alpha_0^{1/(2-q)})+o(\tau^{{2-q}}).
	\end{eqnarray}
	Put (\ref{r1}), (\ref{r2}), and (\ref{r53}) together, we obtain 
	\begin{eqnarray}\nn
	&&2\epsilon^2\alpha_0\tau^{2-2q}q^2D_0+2\epsilon(1-\epsilon)\alpha_0\tau^{2-2q}q^2C_0\\\nn
	&&-2\epsilon(1-\epsilon)f(\rho)\frac{C_qD_q^2}{2-q}\alpha_0^{(q+1)/(2-q)}\phi\left(C_q\left(\frac{1}{1-\rho^2}\right)^{\frac{q}{2(2-q)}}\alpha_0^{1/(2-q)}\right)\\\nn
	&&-2(1-\epsilon)^2\frac{C_qD_q^2}{2-q}\alpha_0^{(q+1)/(2-q)}\phi(C_q\alpha_0^{1/(2-q)})=0,
	\end{eqnarray}
	where $f(\rho)=(\frac{1}{1-\rho^2})^\frac{3q}{2(2-q)}$. Therefore,  for $0\textless\rho\textless 1$, we have
	\begin{eqnarray}\label{rhol1}
	\frac{\tau^{2-2q}}{\alpha_0^{\frac{2q-1}{2-q}}\phi(C_q\alpha_0^{1/(2-q)})}&=&\frac{(1-\epsilon)^2C_qD_q^2}{\epsilon q^2(2-q)(\epsilon D_0+(1-\epsilon)C_0)};
	\end{eqnarray}
	while for $\rho=0$, we have
	\begin{eqnarray}\label{rhoe0}
	\frac{\tau^{2-2q}}{\alpha_0^{\frac{2q-1}{2-q}}\phi(C_q\alpha_0^{1/(2-q)})}&=&\frac{(1-\epsilon)C_qD_q^2}{\epsilon q^2(2-q)C_0}.
	\end{eqnarray}
	From (\ref{rhol1}) and (\ref{rhoe0}), it is straightforward to show that
	\begin{eqnarray}\nn
	e^{-\frac{C_q^2\alpha_0^{\frac{2}{2-q}}}{2}}\sim\tau^{2-2q}&\text{and}&(4-4q)\log\frac{1}{\tau}\sim C_q^2\alpha_0^{\frac{2}{2-q}}.
	\end{eqnarray}
	Combining (\ref{r10}), (\ref{r20}), and (\ref{r30}) together, we obtain the risk function 
	\begin{eqnarray}\nn
	R(\tau^2,\alpha_0\tau^{2-q})&=&\tau^2\left\{\epsilon+\epsilon\alpha_0^2\tau^{2-q}q^2\left(\epsilon D_0+(1-\epsilon)C_0\right)\right\}\\\nn
	&=&\tau^2\left\{\epsilon+\epsilon\left(\frac{4-4q}{C_q^2}\log\frac{1}{\tau}\right)^{2-q}\tau^{2-q}q^2\left(\epsilon D_0+(1-\epsilon)C_0\right)\right\}.
	\end{eqnarray}
	From the high order analysis for the fixed point equation 
	\begin{eqnarray}\nn
	\tau^2=\sigma_w^2+\frac{R(\tau^2,\alpha_0\tau^{2-q})}{\delta},
	\end{eqnarray}
	we obtain the leading order term $\tau^2_0=\frac{\delta}{\delta-\epsilon}\sigma_w^2$ and the next-to-leading order term
	\begin{eqnarray}\nn
	\tau^2_1&=&\frac{q^2\epsilon\tau_0^{4-2q}}{\delta-\epsilon}\left(\frac{4-4q}{C_q^2}\log\frac{1}{\tau_0}\right)^{2-q}\left(\epsilon D_0+(1-\epsilon)C_0\right)\\\nn
	&=&\frac{q^2\epsilon\delta^{2-q}\sigma_w^{4-2q}}{(\delta-\epsilon)^{3-q}}\left(\frac{4-4q}{C_q^2}\log\frac{1}{\sigma_w}\right)^{2-q}\left(\epsilon D_0+(1-\epsilon)C_0\right).
	\end{eqnarray}
	Hence we obtain (3.2).
	
	\begin{figure}[hbtp]
		\begin{center}
			\epsfig{file=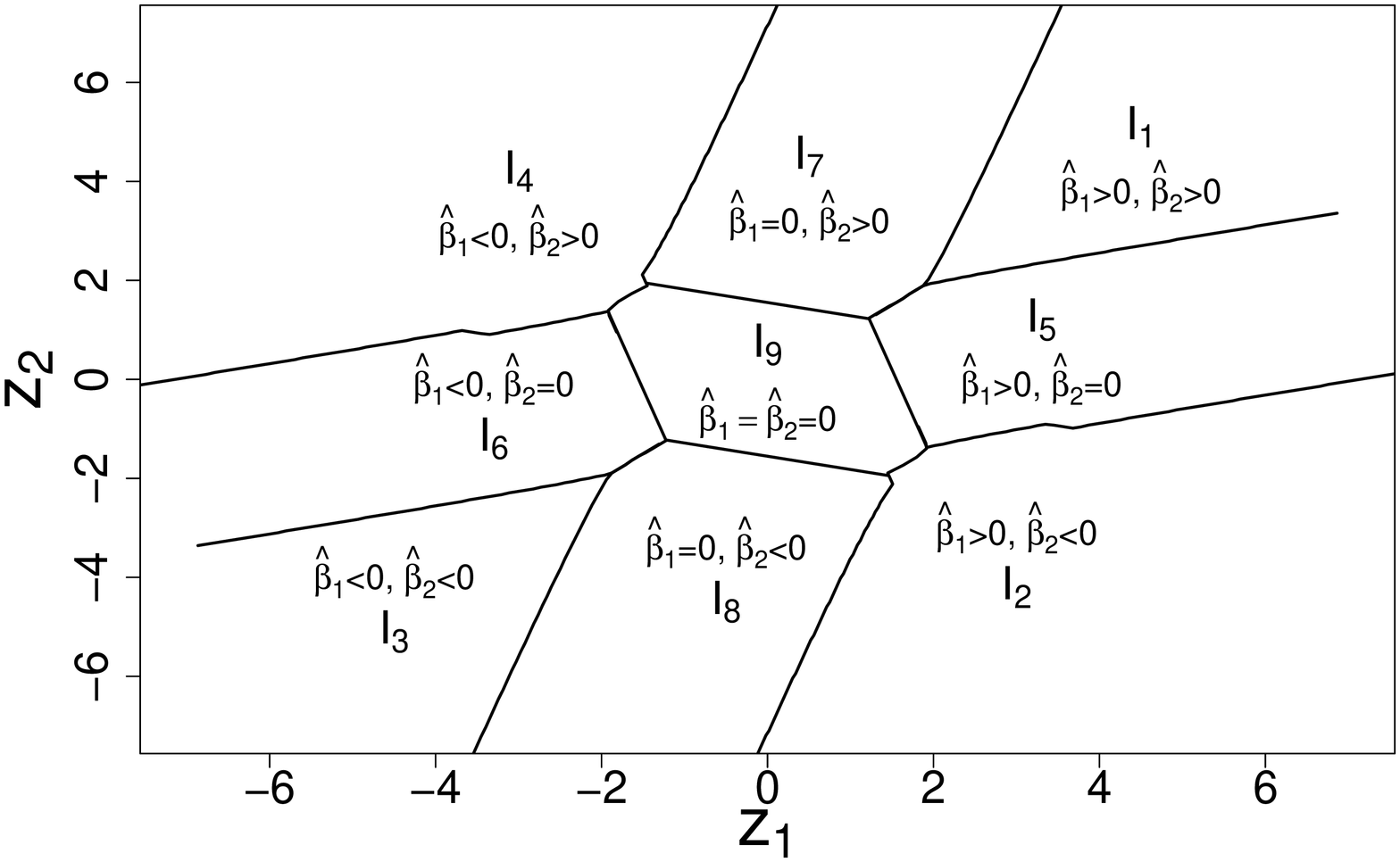,width=8.1cm,angle=-90}
		\end{center}
		\caption{Illustration of the solution for equation (\ref{lfun0}) in two dimensional space. Here $\rho=0.5$ and $\alpha=1$.}
		\label{figure2}
	\end{figure}
\end{proof}

\subsection{Proof of Proposition \ref{prop5}}\label{prop5p}
\begin{proof}
	Denote $U=\bbeta_0/\tau$. For fixed $U$, denote $S=\{j|U_{j}\ne 0\}$. From (\ref{rqd}), we obtain
	\begin{eqnarray}\nn
	R_q(\tau^2,\alpha_0)	&\ge&\frac{1}{p}E[\{(\bSigma^{1/2}\bz)_S-q\alpha_0|\hat{\bbeta}_S|^{q-1}\text{sign}(\hat{\bbeta}_S)\}^T\bSigma_{SS}^{-1}\{(\bSigma^{1/2}\bz)_S-q\alpha_0|\hat{\bbeta}_S|^{q-1}\text{sign}(\hat{\bbeta}_S)\}
	\end{eqnarray} 
	which goes to infinity as $\tau\rightarrow 0$ for fixed $\alpha_0\textgreater 0$ since $\hat{\bbeta}_S\rightarrow\frac{\bbeta_{0,S}}{\tau}$ and $q\textgreater 1$. Therefore, $\alpha_{0\star}(\tau)=0$ and we obtain $R_q(\tau^2,\alpha_\star(\tau))\xrightarrow{\tau\rightarrow 0}1$. From (\ref{fixed}), we have $\tau^2=\sigma_w^2+\frac{\tau^2R_q(\tau^2,\alpha_{0\star}(\tau))}{\delta}$ which implies that $\tau^2(\delta-R_q(\tau^2,\alpha_{0\star}(\tau)))=\delta\sigma_w^2$. If $\delta\textgreater 1$, we have $\tau\xrightarrow{\sigma_w\rightarrow 0}0$ since $R_q(\tau^2,\alpha_{0\star}(\tau))\xrightarrow{\tau\rightarrow 0}1\textless\delta$. On the other hand, if $\delta\textless 1$, we have $\tau\textgreater 0$ as $\sigma_w\rightarrow 0$ because $R_q(\tau^2,\alpha_{0\star}(\tau))\xrightarrow{\tau\rightarrow 0}1$ and $R_q(\tau^2,\alpha_\star(\tau))\xrightarrow{\tau\rightarrow\infty}0$ which leads to $\tau\xrightarrow{\sigma_w\rightarrow 0}\tau_\star\textgreater 0$ such that $R_q(\tau_\star,\alpha_{0\star}(\tau_\star))=\delta$. 
\end{proof}

\subsection{Proof of Proposition \ref{prop6}}
\begin{proof}
	For $\lambda=\alpha_0\tau^{2-q}$, the the solution of $\veta_{q}$ in (\ref{veta}) in the case $1\textless q\le 2$ is the minimizer of loss 
	\begin{eqnarray}\nn
	{\cal L}&=&\frac{1}{2}\|\bbeta-\bbeta_0-\tau\Sigma^{-1/2}\bz\|_\Sigma^2+\alpha_0\tau^{2-q}\|\bbeta\|_q^q.
	\end{eqnarray} 
	It is the solution of the following equation
	\begin{eqnarray}\label{betah0}
	\bSigma(\hat{\bbeta}-\bbeta_0)-\tau\bSigma^{1/2}\bz+q\alpha_0\tau^{2-q}|\hat{\bbeta}|^{q-1}\text{sign}(\hat{\bbeta})=0
	\end{eqnarray}
	which can also be written as 
	\begin{eqnarray}\label{beta0}
	\hat{\bbeta}-\bbeta_0=\tau\bSigma^{-1}\{\bSigma^{1/2}\bz-q\alpha_0\tau^{1-q}|\hat{\bbeta}|^{q-1}\text{sign}(\hat{\bbeta})\},
	\end{eqnarray} 
	or
	\begin{eqnarray}\label{betas}
	\hat{\bbeta}=\bbeta_0+\tau\bSigma^{-1/2}\bz-q\alpha_0\tau^{2-q}\bSigma^{-1}|\hat{\bbeta}|^{q-1}\text{sign}(\hat{\bbeta}).
	\end{eqnarray} 
	Taking derivative of (\ref{betas}) over $\bz$, we have
	\begin{eqnarray}\nn
	\frac{\partial\hat{\bbeta}}{\partial\bz}&=&\tau\bSigma^{-1/2}-q(q-1)\alpha_0\tau^{2-q}\bSigma^{-1}\text{diag}[|\hat{\bbeta}|^{q-2}\text{sign}(\hat{\bbeta})]\frac{\partial\hat{\bbeta}}{\partial\bz}.
	\end{eqnarray} 
	Thus
	\begin{eqnarray}\nn
	\frac{\partial\hat{\bbeta}}{\partial\bz}&=&\tau\left\{\bI+q(q-1)\alpha_0\tau^{2-q}\bSigma^{-1}\text{diag}[|\hat{\bbeta}|^{q-2}\text{sign}(\hat{\bbeta})]\right\}^{-1}\bSigma^{-1/2}.
	\end{eqnarray} 
	Taking derivative of (\ref{betas}) over $\alpha_0$, we have
	\begin{eqnarray}\nn
	\frac{\partial\hat{\bbeta}}{\partial\alpha_0}&=&-q\tau^{2-q}\bSigma^{-1}|\hat{\bbeta}|^{q-1}\text{sign}(\hat{\bbeta})-q(q-1)\alpha_0\tau^{2-q}\bSigma^{-1}|\hat{\bbeta}|^{q-2}\text{sign}(\hat{\bbeta})\frac{\partial\hat{\bbeta}}{\partial\alpha_0}.
	\end{eqnarray} 
	Thus
	\begin{eqnarray}\label{dabeta}
	\frac{\partial\hat{\bbeta}}{\partial\alpha_0}&=&-\left\{\bI+q(q-1)\alpha_0\tau^{2-q}\bSigma^{-1}\text{diag}[|\hat{\bbeta}|^{q-2}\text{sign}(\hat{\bbeta})]\right\}^{-1}q\tau^{2-q}\bSigma^{-1}|\hat{\bbeta}|^{q-1}\text{sign}(\hat{\bbeta}).
	\end{eqnarray}
	Using (\ref{betah0}), we can derive the risk function as
	\begin{eqnarray}\nn
	R(\tau^2,\alpha_0)&=&\frac{1}{p}E\|\hat{\bbeta}-\bbeta_0\|^2_\Sigma\\\nn
	&=&\frac{\tau^2}{p}E\left\{\|\bz\|^2-2q\alpha_0\tau^{1-q}\bz^T\bSigma^{-1/2}|\hat{\bbeta}|^{q-1}\text{sign}(\hat{\bbeta})\right.\\\label{rtot}
	&&\left.+q^2\alpha_0^2\tau^{2-2q}(|\hat{\bbeta}|^{q-1}\text{sign}(\hat{\bbeta}))^T\bSigma^{-1}|\hat{\bbeta}|^{q-1}\text{sign}(\hat{\bbeta})\right\}.
	\end{eqnarray} 
	Then we can conclude that as $\tau\rightarrow 0$, the optimal $\alpha_0^\star\rightarrow 0$ because $\bbeta\rightarrow\bbeta_0$ and $2-2q\textless 0$. Therefore in order for the third term to be finite, we need $\alpha_0^\star\rightarrow 0$. 
	
	In order to obtain the second dominant term, we need to take the first order derivative of $R(\tau^2,\alpha_0)$ over $\alpha_0$ and assume it equal to 0. For $\bSigma=\left(\begin{array}{cc}1&\rho\\\rho&1\end{array}\right)$, the equation can be decomposed into 2-dimensional block. Denote $\rho_1=\frac{1}{2}(\sqrt{1+\rho}+\sqrt{1-\rho})$ and $\rho_2=\frac{1}{2}(\sqrt{1+\rho}-\sqrt{1-\rho})$, we can obtain
	\begin{eqnarray}\nn
	R(\tau^2,\alpha_0)&=&\frac{1}{2}E\|\hat{\bbeta}-\bbeta_0\|^2_\Sigma\\\nn
	&=&\frac{1}{2}E\{(\hat{\beta}_1-\beta_{0,1})^2+(\hat{\beta}_2-\beta_{0,2})^2+2\rho(\hat{\beta}_1-\beta_{0,1})(\hat{\beta}_2-\beta_{0,2})\},
	\end{eqnarray} 
	and the two explicit equations for $\hat{\bbeta}_1$ and $\hat{\bbeta}_2$ 
	\begin{eqnarray}\nn
	\left\{\begin{array}{ccc}
	\hat{\beta}_1-\beta_{0,1}+\rho(\hat{\beta}_2-\beta_{0,2})-\tau\xi_1+\alpha_0\tau^{2-q}q|\hat{\beta}_1|^{q-1}\text{sign}(\hat{\beta}_1)&=&0,\\
	\rho(\hat{\beta}_1-\beta_{0,1})+\hat{\beta}_2-\beta_{0,2}-\tau\xi_2+\alpha_0\tau^{2-q}q|\hat{\beta}_2|^{q-1}\text{sign}(\hat{\beta}_2)&=&0,
	\end{array}\right.
	\end{eqnarray} 
	where $\xi_1=\rho_1z_1+\rho_2z_2$, $\xi_2=\rho_2z_1+\rho_1z_2$. 
	
	We consider three different cases according to the nonzero components of $\beta_{0,1}$ and $\beta_{0,2}$. In the first case: $\beta_{0,1}\ne 0$ and $\beta_{0,2}\ne 0$. Using (\ref{rtot}), we obtain its contribution to the derivative of $R(\tau^2,\alpha_0)$ as
	\begin{eqnarray}\nn
	R_1^\prime(\tau^2,\alpha_0)&=&\epsilon^2\tau^2E\left\{-2q\tau^{1-q}\bz^T\bSigma^{-1/2}|\bbeta|^{q-1}\text{sign}(\bbeta)-2q(q-1)\alpha_0\tau^{1-q}\bz^T\bSigma^{-1/2}|\bbeta|^{q-2}\text{sign}(\bbeta)\frac{\partial\bbeta}{\partial\alpha_0}\right.\\\nn
	&&+2q^2\alpha_0\tau^{2-2q}(|\bbeta|^{q-1}\text{sign}(\bbeta))^T\bSigma^{-1}|\bbeta|^{q-1}\text{sign}(\bbeta)\\\label{r1d}
	&&\left.+2q^2(q-1)\alpha_0^2\tau^{2-2q}(|\bbeta|^{q-1}\text{sign}(\bbeta))^T\bSigma^{-1}|\bbeta|^{q-2}\text{sign}(\bbeta)\frac{\partial\bbeta}{\partial\alpha_0}\right\}.
	\end{eqnarray} 
	Since both components of $\bbeta_0$ are nonzero, as $\tau\rightarrow 0$, (\ref{betas}) can be written as
	\begin{eqnarray}\nn
	\hat{\bbeta}=\bbeta_0+\tau\bSigma^{-1/2}\bz-p\alpha_0\tau^{2-q}\bSigma^{-1}|\bbeta_0|^{q-1}\text{sign}(\bbeta_0)+o_p(\tau^q). 
	\end{eqnarray} 
	Thus, using (\ref{dabeta}), we obtain the two leading terms of (\ref{r1d}) as
	\begin{eqnarray}\nn
	R_1^\prime(\tau^2,\alpha_0)&=&\frac{1}{2}\epsilon^2\tau^2E\{-2q(q-1)\tau^{2-q}Tr[\bSigma^{-1}\text{diag}\{|\bbeta_0|^{q-2}\text{sign}(\bbeta_0)\}]\\\label{r1d}
	&&+2q^2\alpha_0\tau^{2-2q}(|\bbeta_0|^{q-1}\text{sign}(\bbeta_0))^T\bSigma^{-1}|\bbeta_0|^{q-1}\text{sign}(\bbeta_0)\}.
	\end{eqnarray} 
	In the second case, $\beta_{0,1}\ne 0$ and $\beta_{0,2}=0$ (or $\beta_{0,1}=0$ and $\beta_{0,2}\ne 0$). We have
	\begin{eqnarray}\nn
	R_2(\tau^2,\alpha_0)&=&\frac{1}{2}\epsilon(1-\epsilon)E\{(\hat{\beta}_1-\beta_{0,1})^2+\hat{\beta}_2^2+2\rho(\hat{\beta}_1-\beta_{0,1})\hat{\beta}_2\}\\\nn
	&=&\frac{1}{2}\epsilon(1-\epsilon)E\{\hat{\beta}_2^2+[\alpha_0\tau^{2-q}q|\hat{\beta}_1|^{q-1}\text{sign}(\hat{\beta}_1)-\tau\xi_1-\rho\hat{\beta}_2][\alpha_0\tau^{2-q}q|\hat{\beta}_1|^{q-1}\text{sign}(\hat{\beta}_1)-\tau\xi_1+\rho\hat{\beta}_2]\}\\\nn
	&=&\frac{1}{2}\epsilon(1-\epsilon)E\{(1-\rho^2)\hat{\beta}_2^2+[\alpha_0\tau^{2-q}q|\hat{\beta}_1|^{q-1}\text{sign}(\hat{\beta}_1)-\tau\xi_1]^2\}\\\label{r2}
	&=&\frac{1}{2}\epsilon(1-\epsilon)E\{(1-\rho^2)\hat{\beta}_2^2+\alpha_0^2\tau^{4-2q}q^2|\hat{\beta}_1|^{2q-2}+\tau^2\xi_1^2-2\alpha_0\tau^{3-q}q\xi_1|\hat{\beta}_1|^{q-1}\text{sign}(\hat{\beta}_1)\}.
	\end{eqnarray} 
	Taking derivative over $\alpha_0$, using (\ref{dabeta}), we obtain the dominant terms as
	\begin{eqnarray}\nn
	R_2^\prime(\tau^2,\alpha_0)&=&\epsilon(1-\epsilon)E\{(1-\rho^2)\hat{\beta}_2\frac{\partial\hat{\beta}_2}{\partial\alpha_0}+q^2\alpha_0\tau^{4-2q}|\beta_{0,1}|^{2q-2}\}\\\nn
	&=&-\epsilon(1-\epsilon)q\tau^{2-q}E\{|\hat{\beta}_2|^{q}-\rho\hat{\beta}_2|\beta_{0,1}|^{q-1}\text{sign}(\beta_{0,1})-q\alpha_0\tau^{2-q}|\beta_{0,1}|^{2q-2}\}\\\nn
	&=&-\epsilon(1-\epsilon)q\tau^{2-q}E\{|\hat{\beta}_2|^{q}-\frac{q\alpha_0\tau^{2-q}\rho^2|\beta_{0,1}|^{2q-2}}{1-\rho^2}-q\alpha_0\tau^{2-q}|\beta_{0,1}|^{2q-2}\}\\\label{r2d}
	&=&-\epsilon(1-\epsilon)q\tau^{2}E\left\{\frac{|z|^q}{(1-\rho^2)^{q/2}}-\frac{q\alpha_0\tau^{2-2q}|\beta_{0,1}|^{2q-2}}{1-\rho^2}\right\}.
	\end{eqnarray} 
	In the third case, $\beta_{0,1}=0$ and $\beta_{0,2}=0$. Using (\ref{betas}), we obtain
	\begin{eqnarray}\nn
	\hat{\bbeta}=\tau\bSigma^{-1/2}\bz-q\alpha_0\tau^{2-q}\bSigma^{-1}|\hat{\bbeta}|^{q-1}\text{sign}(\hat{\bbeta}).
	\end{eqnarray} 
	Its contribution to risk function is
	\begin{eqnarray}\label{r3}
	R_3(\tau^2,\alpha_0)&=&\frac{1}{2}(1-\epsilon)^2E\|\hat{\bbeta}\|^2_{\bSigma}.
	\end{eqnarray} 
	Taking derivative over $\alpha_0$, we obtain the dominant term as
	\begin{eqnarray}\label{r3d}
	R_3^\prime(\tau^2,\alpha_0)=(1-\epsilon)^2E\hat{\bbeta}^T\bSigma\frac{\partial\hat{\bbeta}}{\partial\alpha_0}=-(1-\epsilon)^2q\tau^{2-q}E\langle|\hat{\bbeta}|^q\rangle
	=-\frac{2(1-\epsilon)^2q\tau^2}{(1-\rho^2)^{q/2}}E\{|z|^q\}.
	\end{eqnarray} 
	Combining (\ref{r1d}), (\ref{r2d}), and (\ref{r3d}) together, we obtain
	\begin{eqnarray}\nn
	&&R_1^\prime(\tau^2,\alpha_0)+2R_2^\prime(\tau^2,\alpha_0)+R_3^\prime(\tau^2,\alpha_0)\\\nn
	&=&-\frac{2(1-\epsilon)q\tau^2}{(1-\rho^2)^{q/2}}E\{|z|^{q}\}+\frac{2\epsilon(1-\epsilon)q^2\tau^{2}\alpha_0\tau^{2-2q}E\{|\beta_{0}|^{2q-2}\}}{2(1-\rho^2)}\\\nn
	&&+\frac{2\epsilon^2q^2\tau^2\alpha_0\tau^{2-2q}}{1-\rho^2}\left[E\left\{|\beta_0|^{2q-2}\right\}-\rho(E|\beta_0|^{q-1}\text{sign}(\beta_0))^2\right].
	\end{eqnarray} 
	Therefore, the optimal tuning $\alpha_{0\star}$ satisfies
	\begin{eqnarray}\label{alphatau}
	\alpha_{0\star}\tau^{2-2q}=\frac{(1-\epsilon)E(|Z|^{q})(1-\rho^2)^{1-q/2}}{\epsilon q\{E(|\beta_{0}|^{2q-2})-\epsilon\rho[E|\beta_0|^{q-1}\text{sign}(\beta_0)]^2\}},
	\end{eqnarray} 
	where $Z\sim N(0,1)$ is independent of $\beta_0$. Now we can obtain the dominant terms in the corresponding $R$ functions. From (\ref{rtot}), we obtain
	\begin{eqnarray}\nn
	R_1(\tau^2,\alpha_{0\star})&=&\tau^2\epsilon^2E\left\{1+\frac{1}{2}q^2\alpha_{0\star}^2\tau^{2-2q}(|\bbeta_0|^{q-1}\text{sign}(\bbeta_0))^T\bSigma^{-1}|\bbeta_0|^{q-1}\text{sign}(\bbeta_0)\right\}.
	\end{eqnarray} 
	From (\ref{r2}), we obtain
	\begin{eqnarray}\nn
	R_2(\tau^2,\alpha_{0\star})&=&2\epsilon(1-\epsilon)\tau^2E\left\{1-q\alpha_{0\star}\tau^{-q}|\hat{\beta}_2|^{q}+\frac{q^2\alpha_{0\star}^2\tau^{2-2q}|\beta_{0}|^{2q-2}}{1-\rho^2}\right\}.
	\end{eqnarray} 
	From (\ref{r3}), we obtain
	\begin{eqnarray}\nn
	R_3(\tau^2,\alpha_{0\star})&=&(1-\epsilon)^2\tau^2E\left\{\frac{1}{2}\|\bz\|^2-q\alpha_{0\star}\tau^{1-q}\bz^T\bSigma^{-1/2}|\hat{\bbeta}|^{q-1}\text{sign}(\hat{\bbeta})\right\}\\\nn
	&=&(1-\epsilon)^2\tau^2E\left\{\frac{1}{2}\|\bz\|^2-q\alpha_{0\star}\tau^{-q}\langle|\hat{\bbeta}|^{q}\rangle\right\}\\\nn
	&=&(1-\epsilon)^2\tau^2\left(1-\frac{q\alpha_{0\star}}{(1-\rho^2)^{q/2}}E\left\{\langle|\bz|^{q}\rangle\right\}\right).
	\end{eqnarray} 
	Putting them together, we have
	\begin{eqnarray}\label{finalrisk}
	R(\tau^2,\alpha_{0\star})&=&\tau^2\left(1+\frac{q^2\alpha_{0\star}^2\tau^{2-2q}}{1-\rho^2}(E|\beta_{0}|^{2q-2}-\epsilon\rho[E|\beta_0|^{q-1}\text{sign}(\beta_0)]^2)\right).
	\end{eqnarray} 
	Substituting (\ref{alphatau}) into (\ref{finalrisk}), we have
	\begin{eqnarray}\nn
	R(\tau^2,\alpha_{0\star})&=&\tau^2+C\tau^{2q},
	\end{eqnarray} 
	where
	\begin{eqnarray}\nn
	C&=&\frac{(1-\epsilon)^2}{\epsilon}\frac{[E(|Z|^{q})]^2(1-\rho^2)^{1-q}\tau^{2q}}{\{E(|\beta_{0}|^{2q-2})-\epsilon\rho[E|\beta_0|^{q-1}\text{sign}(\beta_0)]^2\}}.
	\end{eqnarray} 
	Therefore, the fixed point equation is
	\begin{eqnarray}\nn
	\tau^2=\sigma_w^2+\frac{\tau^2+C\tau^{2q}}{\delta}.
	\end{eqnarray} 
	The first and second dominant terms are
	\begin{eqnarray}\nn
	\tau^2=\frac{\delta}{\delta-1}\sigma_w^2-\frac{C\tau^{2q}}{\delta-1}.
	\end{eqnarray} 
	We obtain the prediction risk
	\begin{eqnarray}\nn
	Risk&=&\delta(\tau^2-\sigma_w^2)=\frac{\delta\sigma_w^2}{\delta-1}-\frac{C\delta\tau^{2q}}{\delta-1}\\\nn
	&=&\frac{\delta\sigma_w^2}{\delta-1}-\frac{\delta^{q+1}}{(\delta-1)^{q+1}}\frac{(1-\epsilon)^2}{\epsilon}\frac{[E(|Z|^{q})]^2(1-\rho^2)^{1-q}\sigma_w^{2q}}{\{E(|\beta_{0}|^{2q-2})-\epsilon\rho[E|\beta_0|^{q-1}\text{sign}(\beta_0)]^2\}}.
	\end{eqnarray} 
\end{proof}

\bibliographystyle{chicago} \bibliography{biblist}

\begin{thebibliography}{}

\bibitem[\protect\citeauthoryear{Amelunxen, Lotz, McCoy, and Tropp}{Amelunxen
  et~al.}{2013}]{Amelunxen2013LivingOT}
Amelunxen, D., M.~Lotz, M.~B. McCoy, and J.~A. Tropp (2013).
\newblock Living on the edge: phase transitions in convex programs with random
  data.
\newblock {\em Information and Inference: A Journal of the IMA\/}~{\em 3},
  224--294.

\bibitem[\protect\citeauthoryear{Baddeley}{Baddeley}{1977}]{baddeley_1977}
Baddeley, A. (1977).
\newblock Integrals on a moving manifold and geometrical probability.
\newblock {\em Advances in Applied Probability\/}~{\em 9\/}(3), 588–603.

\bibitem[\protect\citeauthoryear{Celentano, Montanari, and Wei}{Celentano
  et~al.}{2020}]{wei}
Celentano, M., A.~Montanari, and Y.~Wei (2020).
\newblock The lasso with general gaussian designs with applications to
  hypothesis testing.
\newblock {\em arXiv\/}, 10.48550/ARXIV.2007.13716.

\bibitem[\protect\citeauthoryear{Chartrand and Staneva}{Chartrand and
  Staneva}{2008}]{Chartrand_2008}
Chartrand, R. and V.~Staneva (2008, may).
\newblock Restricted isometry properties and nonconvex compressive sensing.
\newblock {\em Inverse Problems\/}~{\em 24\/}(3), 035020.

\bibitem[\protect\citeauthoryear{Donoho}{Donoho}{2006}]{https://doi.org/10.1002/cpa.20132}
Donoho, D.~L. (2006).
\newblock For most large underdetermined systems of linear equations the
  minimal $\ell_1$-norm solution is also the sparsest solution.
\newblock {\em Communications on Pure and Applied Mathematics\/}~{\em 59\/}(6),
  797--829.

\bibitem[\protect\citeauthoryear{Donoho, Maleki, and Montanari}{Donoho
  et~al.}{2011}]{dmm}
Donoho, D.~L., A.~Maleki, and A.~Montanari (2011, Oct).
\newblock The noise-sensitivity phase transition in compressed sensing.
\newblock {\em IEEE Transactions on Information Theory\/}~{\em 57\/}(10),
  6920--6941.

\bibitem[\protect\citeauthoryear{Donoho and Tanner}{Donoho and
  Tanner}{2005}]{Donoho9446}
Donoho, D.~L. and J.~Tanner (2005).
\newblock Sparse nonnegative solution of underdetermined linear equations by
  linear programming.
\newblock {\em Proceedings of the National Academy of Sciences\/}~{\em
  102\/}(27), 9446--9451.

\bibitem[\protect\citeauthoryear{Fu and Knight}{Fu and
  Knight}{2000}]{10.1214/aos/1015957397}
Fu, W. and K.~Knight (2000).
\newblock {Asymptotics for lasso-type estimators}.
\newblock {\em The Annals of Statistics\/}~{\em 28\/}(5), 1356 -- 1378.

\bibitem[\protect\citeauthoryear{Huang}{Huang}{2021}]{lasso}
Huang, H. (2021).
\newblock Lasso risk and phase transition under dependence.
\newblock {\em arXiv\/}, 10.48550/ARXIV.2103.16035.

\bibitem[\protect\citeauthoryear{Javanmard and Montanari}{Javanmard and
  Montanari}{2014}]{Montanari}
Javanmard, A. and A.~Montanari (2014, Oct).
\newblock Hypothesis testing in high-dimensional regression under the gaussian
  random design model: Asymptotic theory.
\newblock {\em IEEE Transactions on Information Theory\/}~{\em 60\/}(10),
  6522--6554.

\bibitem[\protect\citeauthoryear{Kabashima, Wadayama, and Tanaka}{Kabashima
  et~al.}{2009}]{Kabashima_2009}
Kabashima, Y., T.~Wadayama, and T.~Tanaka (2009, sep).
\newblock A typical reconstruction limit for compressed sensing based on
  $\ell_p$-norm minimization.
\newblock {\em Journal of Statistical Mechanics: Theory and Experiment\/}~{\em
  2009\/}(09), L09003.

\bibitem[\protect\citeauthoryear{Ma, Xu, and Maleki}{Ma et~al.}{2019}]{8618402}
Ma, J., J.~Xu, and A.~Maleki (2019).
\newblock Optimization-based amp for phase retrieval: The impact of
  initialization and $\ell_{2}$ regularization.
\newblock {\em IEEE Transactions on Information Theory\/}~{\em 65\/}(6),
  3600--3629.

\bibitem[\protect\citeauthoryear{Meinshausen and B{\"u}hlmann}{Meinshausen and
  B{\"u}hlmann}{2006}]{10.1214/009053606000000281}
Meinshausen, N. and P.~B{\"u}hlmann (2006).
\newblock {High-dimensional graphs and variable selection with the Lasso}.
\newblock {\em The Annals of Statistics\/}~{\em 34\/}(3), 1436 -- 1462.

\bibitem[\protect\citeauthoryear{M{\'e}zard and Montanari}{M{\'e}zard and
  Montanari}{2009}]{mezard2009information}
M{\'e}zard, M. and A.~Montanari (2009).
\newblock {\em Information, Physics, and Computation}.
\newblock Oxford Graduate Texts. OUP Oxford.

\bibitem[\protect\citeauthoryear{Mezard, Parisi, and Virasoro}{Mezard
  et~al.}{1987}]{mezard1987spin}
Mezard, M., G.~Parisi, and M.~Virasoro (1987).
\newblock {\em Spin Glass Theory and Beyond: An Introduction to the Replica
  Method and Its Applications}.
\newblock World Scientific Lecture Notes in Physics. World Scientific.

\bibitem[\protect\citeauthoryear{Rangan, Goyal, and Fletcher}{Rangan
  et~al.}{2009}]{NIPS2009_bdb106a0}
Rangan, S., V.~Goyal, and A.~K. Fletcher (2009).
\newblock Asymptotic analysis of map estimation via the replica method and
  compressed sensing.
\newblock In Y.~Bengio, D.~Schuurmans, J.~Lafferty, C.~Williams, and A.~Culotta
  (Eds.), {\em Advances in Neural Information Processing Systems}, Volume~22.
  Curran Associates, Inc.

\bibitem[\protect\citeauthoryear{Saab, Chartrand, and Yilmaz}{Saab
  et~al.}{2008}]{4518502}
Saab, R., R.~Chartrand, and O.~Yilmaz (2008).
\newblock Stable sparse approximations via nonconvex optimization.
\newblock In {\em 2008 IEEE International Conference on Acoustics, Speech and
  Signal Processing}, pp.\  3885--3888.

\bibitem[\protect\citeauthoryear{Saab and Özgür Yılmaz}{Saab and Özgür
  Yılmaz}{2010}]{SAAB201030}
Saab, R. and Özgür Yılmaz (2010).
\newblock Sparse recovery by non-convex optimization – instance optimality.
\newblock {\em Applied and Computational Harmonic Analysis\/}~{\em 29\/}(1),
  30--48.

\bibitem[\protect\citeauthoryear{Stojnic}{Stojnic}{2013}]{https://doi.org/10.48550/arxiv.1306.3976}
Stojnic, M. (2013).
\newblock Lifting $\ell_q$-optimization thresholds.
\newblock {\em arXiv\/}, https://arxiv.org/abs/1306.3976.

\bibitem[\protect\citeauthoryear{Thrampoulidis, Abbasi, and
  Hassibi}{Thrampoulidis et~al.}{2018}]{8365826}
Thrampoulidis, C., E.~Abbasi, and B.~Hassibi (2018).
\newblock Precise error analysis of regularized $m$ -estimators in high
  dimensions.
\newblock {\em IEEE Transactions on Information Theory\/}~{\em 64\/}(8),
  5592--5628.

\bibitem[\protect\citeauthoryear{Wang, Xu, and Tang}{Wang
  et~al.}{2011}]{5893944}
Wang, M., W.~Xu, and A.~Tang (2011).
\newblock On the performance of sparse recovery via $\ell_p$-minimization $(0
  \leq p \leq 1)$.
\newblock {\em IEEE Transactions on Information Theory\/}~{\em 57\/}(11),
  7255--7278.

\bibitem[\protect\citeauthoryear{Wang, Weng, and Maleki}{Wang
  et~al.}{2020}]{10.1214/19-AOS1906}
Wang, S., H.~Weng, and A.~Maleki (2020).
\newblock {Which bridge estimator is the best for variable selection?}
\newblock {\em The Annals of Statistics\/}~{\em 48\/}(5), 2791 -- 2823.

\bibitem[\protect\citeauthoryear{Wang, Weng, and Maleki}{Wang
  et~al.}{2021}]{10.1093/imaiai/iaab025}
Wang, S., H.~Weng, and A.~Maleki (2021).
\newblock {Does SLOPE outperform bridge regression?}
\newblock {\em Information and Inference: A Journal of the IMA\/}~{\em
  11\/}(1), 1--54.

\bibitem[\protect\citeauthoryear{Weng and Maleki}{Weng and
  Maleki}{2019}]{10.1093/imaiai/iay024}
Weng, H. and A.~Maleki (2019).
\newblock {Low noise sensitivity analysis of -minimization in oversampled
  systems}.
\newblock {\em Information and Inference: A Journal of the IMA\/}~{\em 9\/}(1),
  113--155.

\bibitem[\protect\citeauthoryear{Weng, Maleki, and Zheng}{Weng
  et~al.}{2018}]{10.1214/17-AOS1651}
Weng, H., A.~Maleki, and L.~Zheng (2018).
\newblock {Overcoming the limitations of phase transition by higher order
  analysis of regularization techniques}.
\newblock {\em The Annals of Statistics\/}~{\em 46\/}(6A), 3099 -- 3129.

\bibitem[\protect\citeauthoryear{Weng, Zheng, Maleki, and Wang}{Weng
  et~al.}{2016}]{7541384}
Weng, H., L.~Zheng, A.~Maleki, and X.~Wang (2016).
\newblock Phase transition and noise sensitivity of $l_p$-minimization for
  $0\le p\le 1$.
\newblock In {\em 2016 IEEE International Symposium on Information Theory
  (ISIT)}, pp.\  675--679.

\bibitem[\protect\citeauthoryear{Zheng, Maleki, Weng, Wang, and Long}{Zheng
  et~al.}{2017}]{7987040}
Zheng, L., A.~Maleki, H.~Weng, X.~Wang, and T.~Long (2017).
\newblock Does $\ell _{p}$ -minimization outperform $\ell _{1}$ -minimization?
\newblock {\em IEEE Transactions on Information Theory\/}~{\em 63\/}(11),
  6896--6935.

\end{thebibliography}
\end{document}